\documentclass[review,11pt]{elsarticle}

\usepackage{lineno,hyperref}
\modulolinenumbers[5]

\usepackage{graphicx}
\usepackage{bm}
\usepackage{color}
\usepackage{anysize}
\usepackage{booktabs}
\usepackage[linesnumbered,ruled]{algorithm2e}
\usepackage{subfig}
\usepackage{xcolor}
\usepackage{amsmath,amsfonts,amsthm}
\usepackage{multirow}
\usepackage{longtable}
\usepackage{float}
\usepackage{url}
\usepackage{hyperref}
\usepackage{lscape}
\usepackage{caption}

\newtheorem{lemma}{Lemma}
\newtheorem{theo}{Theorem}

\journal{European Journal of Operational Research}

\usepackage{pdfpages}
\begin{document}
\linespread{1.5}
\includepdf{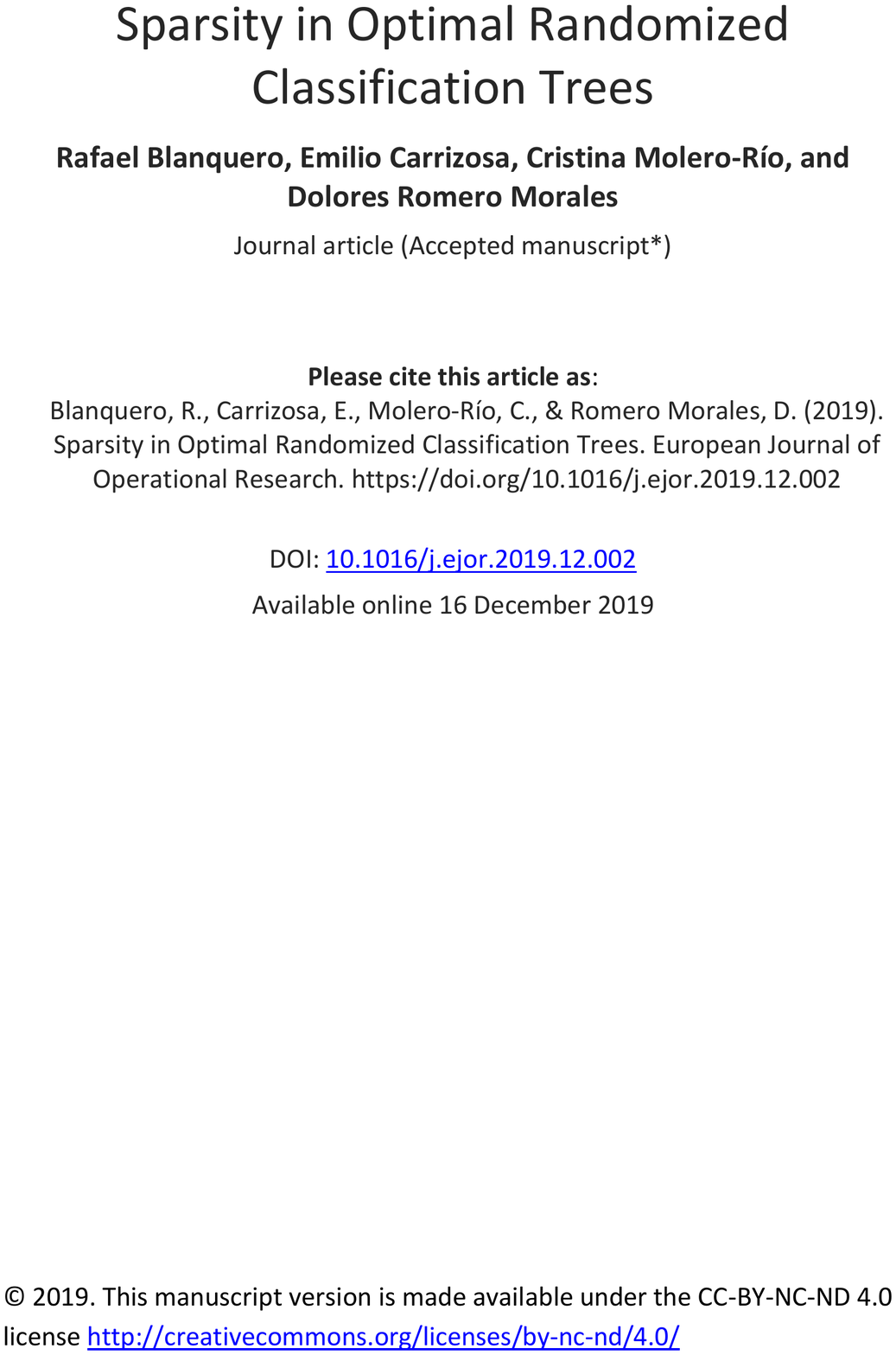}
\begin{frontmatter}

\title{Sparsity in Optimal Randomized Classification Trees}


\author[mymainaddress]{Rafael Blanquero}
\ead{rblanquero@us.es}

\author[mymainaddress]{Emilio Carrizosa}
\ead{ecarrizosa@us.es}

\author[mymainaddress]{Cristina Molero-R\'io\corref{mycorrespondingauthor}}
\cortext[mycorrespondingauthor]{Corresponding author}
\ead{mmolero@us.es}

\author[mysecondaryaddress]{Dolores Romero Morales}
\ead{drm.eco@cbs.dk}

\address[mymainaddress]{Instituto de Matem\'aticas de la Universidad de Sevilla (IMUS), Seville, Spain}
\address[mysecondaryaddress]{Copenhagen Business School, Frederiksberg, Denmark}

\begin{abstract}
Decision trees are popular Classification and Regression tools and, when small-sized, easy to interpret. Traditionally, a greedy approach has been used to build the trees, yielding a very fast training process; however, controlling sparsity (a proxy for interpretability) is challenging. In recent studies, optimal decision trees, where all decisions are optimized simultaneously, have shown a better learning performance, especially when oblique cuts are implemented. In this paper, we propose a continuous optimization approach to build sparse optimal classification trees, based on oblique cuts, with the aim of using fewer predictor variables in the cuts as well as along the whole tree. Both types of sparsity, namely local and global, are modeled by means of regularizations with polyhedral norms.  The computational experience reported supports the usefulness of our methodology. In all our data sets, local and global sparsity can be improved without harming classification accuracy. Unlike greedy approaches, our ability to easily trade in some of our classification accuracy for a gain in global sparsity is shown.

\end{abstract}

\begin{keyword}
 Data mining \sep Optimal Classification Trees \sep Global and Local Sparsity \sep Nonlinear Programming
\end{keyword}

\end{frontmatter}

\section{Introduction}
\label{Introduction}

Decision trees \cite{yang2017regression} are a popular non-parametric tool for Classification and Regression in Statistics and Machine Learning \cite{hastie2009elements}. Since they are rule-based, when small-sized, they are deemed to be leaders in terms of interpretability \cite{athey2018impact,BaesensMS03,carrizosa2011detecting,freitasACM14,goodman16,jung2017simple,martens2007comprehensible,martin2014interpretable,ridgewayNIJJ13,ustunML16}.

It is well-known that the problem of building optimal decision trees is NP-complete \cite{hyafil1976constructing}.

For this reason, classic decision trees have been traditionally designed using greedy procedures in which at each branch node of the tree, some purity criterion is (locally) optimized. For instance, CARTs \cite{breiman1984classification} employ a greedy and recursive partitioning procedure which is computationally cheap, especially since orthogonal cuts are implemented, i.e., one single predictor variable is involved in each branching rule. These rules are of maximal sparsity at each branching node (excellent local sparsity), making classic decision trees locally easy to interpret. However, when deep, they become to be harder to interpret since many predictor variables are, in general, involved across all branching rules (not so good global sparsity).

Addressing global sparsity is a challenge in decision trees and, to the best of our knowledge, this has not been tackled appropriately in the literature. Standard CARTs or Random Forests (RFs) \cite{biau2016random,breiman2001random,fernandez2014we,genuer2017random} cannot manage it due to the greedy construction of the trees. Nonetheless, some attempts have been made, see \cite{deng2012feature,deng2013gene}. Classic decision trees usually select their orthogonal cuts at each branch node by optimizing an information theory criterion among all possible predictor variables and thresholds. The regularization framework in \cite{deng2012feature} considers a penalty to this criterion for predictor variables that have not appeared yet in the tree. This approach is refined in \cite{deng2013gene}, by also including the importance scores of the predictor variables, obtained in a preprocessing step running a preliminary RF.

The mainstream trend of using a greedy strategy in the construction of decision trees may lead to myopic decisions, which, in turn, may affect the overall learning performance. The major advances in Mathematical Optimization \cite{carrizosa2013supervised,olafsson2008operations,silva2017optimization} have led to different approaches to build decision trees with some overall optimality criterion, called hereafter optimal classification trees. It is worth mentioning recent proposals which grow optimal classification trees of a pre-established depth, both deterministic \cite{bertsimas2017optimal,firat2018constructing,menickelly2016optimal,verwer2017learning,verwer2017auction} and randomized \cite{orct}. The deterministic approaches formulate the problem of building the tree as a mixed-integer linear optimization problem. Such approach is the most natural, since many discrete decisions are to be made when building a decision tree. Although the results of such optimal classification trees are encouraging, the inclusion of integer decision variables makes the computing times explode, giving rise to models trained over a small subsample of the data set \cite{menickelly2016optimal} and, as customary, with a CPU time limit being imposed to the optimization solver. On the other hand, a continuous optimization-based approach to build optimal  randomized classification trees is proposed in \cite{orct}. This is achieved by replacing the yes/no decisions in traditional trees by probabilistic decisions, i.e., instead of deciding at each branch node if an individual goes either to the left or to the right child node in the tree, the probability of going to the left is sought. The numerical results in \cite{orct} illustrate the good performance achieved in very short time. All these optimization-based approaches are flexible enough to address critical issues that the greedy nature of classic decision trees would find it difficult, such as preferences on the classification performance in some class where misclassifying is more damaging \cite{orct,verwer2017learning,verwer2017auction}, or controlling the number of predictor variables used along the tree (local and global sparsity).

Optimal classification trees have been grown with both orthogonal \cite{bertsimas2017optimal,firat2018constructing,menickelly2016optimal} and oblique cuts \cite{bennett1996optimal,bertsimas2017optimal,orct,norouzi2015efficient,verwer2017learning,verwer2017auction}.
Oblique cuts are more flexible than orthogonal ones since a combination of several predictor variables is allowed in the branching. Trees based on oblique cuts lead to similar or even better learning performance than those based on orthogonal cuts, and, at the same time, they exhibit a shallow depth, since several orthogonal cuts may be reduced to one single oblique cut. Apart from the flexibility that we can borrow from them, many integer decision variables associated with orthogonal cuts are not present in the oblique ones, which eases the optimization. Therefore, optimal classification trees based on oblique cuts require a lower training computing time while showing much more promising results in terms of accuracy. However, this comes at the expense of damaging interpretability, since, in principle, all the predictor variables could appear in each branching rule. In this paper, we tackle this issue.

We propose a novel optimized classification tree, based on the methodology in \cite{orct} and, therefore, in oblique cuts, that yields rules/trees that are sparser, and thus enhance interpretability. We model this as a continuous optimization problem.

As in the classic LASSO model \cite{tibshirani2015statistical}, sparsity is sought by means of regularization terms.

 We model local sparsity with the $\ell_1$-norm, and the global sparsity with the $\ell_\infty$-norm. The $\ell_\infty$ reguralization has been applied to other classifiers, for instance, Support Vector Machines \cite{maldonado2017integrated,maldonado2017synchronized,zou2008f}, but the $\ell_1$ is more popular. A novel continuous-based approach for building this sparse optimal randomized classification tree is provided. Theoretical results on the range of the sparsity parameters are shown. Our numerical results, where well-known real data sets are used, illustrate the efectiveness of our methodology: sparsity in optimal classification trees improves without harming learning performance. In addition, our ability to trade in some of our classification accuracy, still being superior to CART, to be comparable to CART in terms of global sparsity is shown.

The remainder of the paper is organized as follows. In Section \ref{Sparsity in Optimal Randomized Classification Trees} we detail the construction of the Sparse Optimal Randomized Classification Tree. Some theoretical properties are given in Section \ref{Theoretical properties}. In Section \ref{Computational experience}, our numerical experience is reported.

Finally, conclusions and possible lines of future research are provided in Section \ref{Conclusions and future research}.

\section{Sparsity in Optimal Randomized Classification Trees}
\label{Sparsity in Optimal Randomized Classification Trees}

\subsection{Introduction}
\label{Introduction}

We assume given a training sample $\left\lbrace \left( \bm{x}_i,y_i \right) \right\rbrace_{1\leq i\leq N}$, where $\bm{x}_i$ represents the $p$-dimensional vector of predictor variables of individual $i$, and $y_i\in\left\lbrace 1,\ldots, K\right\rbrace$ indicates the class membership. Without loss of generality, we assume $\bm{x}_i\in\left[ 0,1\right]^p,\,\,i=1,\ldots,N$.

\begin{figure}[h!]
	\centering
	\includegraphics[scale=0.45]{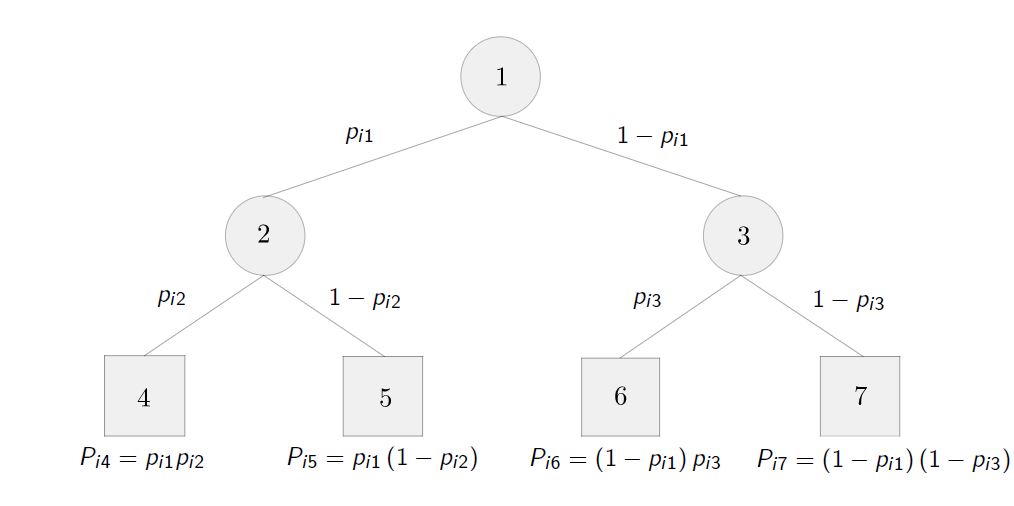}
	\caption{Optimal Randomized Classification Tree of depth $D=2$.}
	\label{Optimal randomized classification tree for depth $D=2$.}
\end{figure}

Sparse Optimal Randomized Classification Trees, addressed in this paper, extend the Optimal Randomized Classification Trees (ORCTs) in \cite{orct}. An ORCT is an optimal binary classification tree of a given depth $D$, obtained by minimizing the expected misclassification cost over the training sample. Figure \ref{Optimal randomized classification tree for depth $D=2$.} 
shows the structure of an ORCT of depth $D=2$. Unlike classic decision trees, oblique cuts, on which more than one predictor variable takes part, are performed. ORCTs are modeled by means of a Non-Linear Continuous Optimization formulation. The usual deterministic yes/no rule at each branch node is replaced by a smoother rule: a probabilistic decision rule at each branch node, induced by a cumulative density function (CDF) $F$, is obtained. Therefore, the movements in ORCTs can be seen as randomized: at a given branch node of an ORCT, a random variable will be generated to indicate by which branch an individual has to continue. Since binary trees are built, the Bernoulli distribution is appropriate, whose probability of success will be determined by the value of this CDF, evaluated over the vector of predictor variables. More precisely, at a given branch node $t$ of the tree, an individual with predictor variables $\bm{x}$ will go either to the left or to the right child nodes with probabilities $F\left(\dfrac{1}{p} \bm{a}_{\cdot t}^T \bm{x}-\mu_t\right)$ and $1-F\left(\dfrac{1}{p}\bm{a}_{\cdot t}^T \bm{x}-\mu_t\right)$, respectively, where $\bm{a}_{\cdot t}$ and $\mu_t$ are decision variables. For further details on the construction of ORCTs, the reader is referred to \cite{orct}. Sparse ORCT, S-ORCT, minimizes the expected misclassification cost over the training sample regularized with two polyhedral norms.

The following notation is needed:

\begin{longtable}{l p{12.25cm}}
\textit{Parameters} & \\
$D$ & depth of the binary tree, \\
$N$ & number of individuals in the training sample,\\
$p$ & number of predictor variables, \\
$K$ & number of classes, \\
$\left\lbrace \left( \bm{x}_i,y_i\right)\right\rbrace_{1\leq i\leq N}$ & training sample, where $\bm{x}_i \in \left[0,1\right]^p$ and $y_i\in \left\lbrace 1,\ldots,K\right\rbrace,$\\
$I_k$ & set of individuals in the training sample belonging to class $k$, $k=1,\ldots,K,$\\
$W_{y_ik}$ & misclassification cost incurred when classifying an individual $i$, whose class is $y_i$, in class $k$, $y_i, i=1,\ldots,N, \,\,k=1,\ldots,K$, \\
$F\left(\cdot\right)$ & univariate continuous CDF centered at $0$, used to define the probabilities for an individual to go to the left or the right child node in the tree. We will assume that $F$ is the CDF of a continuous random variable with density $f$,\\
$\lambda^L\geq 0$ & local sparsity regularization parameter,\\
$\lambda^G\geq 0$ & global sparsity regularization parameter,\\
& \\
\textit{Nodes} &\\

$\tau_B$ & set of branch nodes,\\
$\tau_L$ & set of leaf nodes,\\
$N_L\left( t\right)$ & set of ancestor nodes of leaf node $t$ whose left branch takes part in the path from the root node to leaf node $t$, $t\in\tau_L$,\\
$N_R\left( t\right)$ & set of ancestor nodes of leaf node $t$ whose right branch takes part in the path from the root node to leaf node $t$, $t\in\tau_L$,\\
& \\

\textit{Decision variables} & \\

$a_{jt}\in\left[-1,1\right]$ & coefficient of predictor variable $j$ in the oblique cut at branch node $t\in\tau_B$, with $\bm{a}$ being the $p\times\left\lvert\tau_B\right\rvert$ matrix of these coefficients, $\bm{a}=\left( a_{jt}\right)_{j=1,\ldots,p,\,\,t\in\tau_B}$. The expressions $\bm{a}_{j\cdot}$ and $\bm{a}_{\cdot t}$ will denote the $j$-th row and the $t$-th column of $\bm{a}$, respectively, \\
$\mu_t\in\left[-1,1\right]$ & location parameter at branch node $t\in\tau_B$, $\bm{\mu}$ being the vector that comprises every $\mu_t$, i.e., $\bm{\mu}=\left( \mu_t \right)_{t\in\tau_B}$,\\
$C_{kt}$ & probability  of being assigned to class label $k\in\left\lbrace 1,\ldots, K\right\rbrace$ for an individual at leaf node $t,\,\,t\in\tau_L$, being the $K\times\left\lvert\tau_L\right\rvert$ matrix such that $\bm{C}=\left( C_{kt}\right)_{k=1,\ldots,K,\,\,t\in\tau_L}$.\\

& \\

\textit{Probabilities} & \\

$p_{it}\left( \bm{a}_{\cdot t}, \mu_t\right)$ & probability of individual $i$ going down the left branch at branch node $t$. Its expression is $ p_{it}\left( \bm{a}_{\cdot t}, \mu_t\right)= F\left(\dfrac{1}{p} \bm{a}_{\cdot t}^T \bm{x}_{i}-\mu_t\right), \,\, i =1,\ldots, N, \,\, t\in \tau_B$,\\
$P_{it}\left( \bm{a}, \bm{\mu}\right)$ & probability of individual $i$ falling into leaf node $t$. Its expression is
 $P_{it} \left( \bm{a}, \bm{\mu}\right) = \prod\limits_{t_l\in N_L(t)} p_{it_l}\left( \bm{a}_{\cdot t_l}, \mu_{t_l}\right)\prod\limits_{t_r\in N_R(t)}\left(1- p_{it_r}\left( \bm{a}_{\cdot t_r}, \mu_{t_r}\right)\right),\,\, i=1,\ldots,N,\,\, t\in \tau_L,$ \\

$g\left(\bm{a},\bm{\mu},\bm{C}\right)$& expected misclassification cost over the training sample. Its expression is $ g\left(\bm{a},\bm{\mu},\bm{C}\right) = \dfrac{1}{N}\sum\limits_{i=1}^N \sum\limits_{t\in\tau_L} P_{it}\left( \bm{a}, \bm{\mu}\right) \sum\limits_{k=1}^K W_{y_ik} C_{kt}$.\\
\end{longtable}

\subsection{The formulation}
\label{The formulation}

With these parameters and decision variables, the S-ORCT is formulated as follows:
\begin{align}
\min \hspace*{0.6cm}& g\left(\bm{a},\bm{\mu},\bm{C}\right)  +\lambda^L \sum_{j=1}^p \left\lVert  \bm{a}_{j\cdot}\right\rVert_1 +\lambda^G \sum_{j=1}^p \left\lVert \bm{a}_{j\cdot}\right\rVert_\infty  \label{objetivo}\\
\text{s.t. \hspace*{0.5cm}}
& \sum_{k=1}^K C_{kt} = 1, \,\, t\in \tau_L,\label{res2}\\
     & \sum_{t\in \tau_L} C_{kt} \geq 1, \,\, k=1,\ldots, K,\label{res3} \\
     & a_{jt} \in \left[ -1,1\right],\,\, j=1,\ldots,p,\,\, t\in \tau_B,\label{res4}\\
     & \mu_{t} \in \left[-1,1\right] , \,\, t\in \tau_B,\label{res5}\\
     & C_{kt} \in \left[ 0,1\right], \,\, k=1,\ldots,K,\,\, t\in \tau_L.\label{res6}
\end{align}
In the objective function we have three terms, the first being the expected misclassification cost in the training sample, while the second and the third are regularization terms. The second term addresses local sparsity, since it penalizes the coefficients of the predictor variables used in the cuts along the tree. Instead, the third term controls whether a given predictor variable is ever used across the whole tree, thus addressing global sparsity. The $\ell_\infty$-norm is used as a group penalty function, by forcing the coefficients linked to the same predictor variable to be shrunk simultaneously along all branch nodes. Note that both local and global sparsity are equivalent when dealing with depth $D=1$, as there is a single cut across the whole tree.

In terms of the feasible region, for each leaf node $t\in\tau_L$, $C_{kt}$ represents the probability that an individual at node $t$ is assigned to class $k\in\left\lbrace 1,\ldots,K\right\rbrace$. Constraints \eqref{res2} force that such probabilities sum to $1$, while constraints \eqref{res3} force the sum of the probabilities along all leaf nodes $t\in\tau_B$ assigned to class $k$ to be at least one. 

Theorem \ref{resultadosoluciondeterminista} guarantees the existence of an optimal deterministic solution, i.e., such probabilities $C_{kt}$ will all be in $\left\lbrace 0,1\right\rbrace$, and thus \eqref{res6} can be replaced by
\begin{equation}
\label{siete}
C_{kt}\in\left\lbrace 0,1\right\rbrace,\,\, k=1,\ldots,K,\,\, t\in\tau_L.
\end{equation} Constraints \eqref{res6} and \eqref{siete} will be used interchangeably when needed.

\begin{theo}
\label{resultadosoluciondeterminista}
There exists an optimal solution to \eqref{objetivo}-\eqref{res6} such that $C_{kt}\in\left\lbrace 0,1\right\rbrace,\,\, k=1,\ldots,K,\,\,t\in\tau_L$.
\end{theo}

\noindent\textit{Proof.}

The continuity of the objective function \eqref{objetivo}, defined over a compact set, ensures the existence of an optimal solution of the optimization problem \eqref{objetivo}-\eqref{res6}, by Weierstrass Theorem. Let $\bm{a}^*=\left( a^*_{jt}\right)_{j=1,\ldots,p,\,\,t\in\tau_B},\,\, \bm{\mu}^*=\left( \mu_t^*\right)_{t\in\tau_B},\,\, \bm{C}^*=\left(C^*_{kt}\right)_{k=1,\ldots,K,\,\,t\in\tau_B}$ be an optimal solution. Fixed $\bm{a}^*,\,\, \bm{\mu}^*$, then  $\bm{C}^*$ is optimal to the following problem in the decision variables $C_{kt},\,\, k=1,\ldots,K,\,\,t\in\tau_L$:
\begin{align}
\begin{split}
\nonumber
\min \hspace*{0.6cm}& \dfrac{1}{N} \sum_{i=1}^N \sum_{t\in\tau_L} P_{it}\left( \bm{a}^*, \bm{\mu}^*\right) \sum_{k=1}^K W_{y_ik} C_{kt} +\lambda^L \sum_{j=1}^p \left\lVert \bm{a}^*_{j\cdot}
\right\rVert_1 +\lambda^G \sum_{j=1}^p \left\lVert \bm{a}^*_{j\cdot}\right\rVert_\infty\\
\text{s.t. \hspace*{0.5cm}}
     & \sum_{k=1}^K C_{kt} = 1, \,\, t\in \tau_L,\\
     & \sum_{t\in \tau_L} C_{kt} \geq 1, \,\, k=1,\ldots, K, \\
     & C_{kt} \in \left[ 0,1\right], \,\, k=1,\ldots,K,\,\, t\in \tau_L.
 \end{split}
\end{align}
This is a transportation problem, to which the integrality of an optimal solution is well-known to hold, i.e., there exists $\overline{\bm{C}}=\left(\overline{C}_{kt}\right)_{k=1,\ldots,K,\,\, t\in\tau_L}\in\left\lbrace 0,1\right\rbrace$ for all $k,\,\,t$ such that $\left( \bm{a}^*,\bm{\mu}^*, \overline{\bm{C}} \right)$ is also optimal for \eqref{objetivo}-\eqref{res6}.
\qed

Theorem \ref{resultadosoluciondeterminista} gives a new interpretation of constraints \eqref{res2}-\eqref{res3}: if \eqref{siete} is used instead of \eqref{res6}, when $C_{kt}$ takes the value $1$, then all the individuals at node $t\in\tau_L$ are labelled as $k$; and $0$, otherwise. Constraints \eqref{res2} state that any leaf node $t\in\tau_L$ must be labelled with exactly one class label, and constraints \eqref{res3} state that each class $k$ has at least one node $t$ with such label.

Once the optimization problem is solved, the S-ORCT predicts the class of a new unlabeled observation with predictor vector $\bm{x}$ with a probabilistic rule, namely, we estimate the probability of being in class $k$ as $\sum\limits_{t\in \tau_L} C_{kt}\cdot P_{\bm{x}t}\left( \bm{a},\bm{\mu}\right)$. If a deterministic classification rule is sought, we allocate to the most probable class. Moreover, if prior probabilities $\Pi_{k}\left(\bm{x}\right)$ are given, one can also use the Bayes rule.

ORCTs were also shown to deal effectively with controlling the correct classification rate on different classes. This idea can also be applied to S-ORCTs. Hence, given the classes $k=1,\ldots,K$ to be controlled and their corresponding desired performances $\rho_k$, the expectation of achieving each performance guarantee can be computed with the ORCT parameters, provided that the following set of constrainsts is added to the model:
\begin{equation}
\label{performanceconstraints}
\sum_{i\in I_k} \sum_{t\in\tau_L} P_{it}\left( \bm{a},\bm{\mu}\right)C_{kt} \geq \rho_k \vert I_k \vert,\,\, k=1,\ldots,K.
\end{equation}

With these constraints we have a direct control on the classification performance in each class separately. This is useful when dealing with imbalanced data sets.

\subsection{A smooth reformulation}
\label{A smooth reformulation}

Problem \eqref{objetivo}-\eqref{res6} is non-smooth due to the norms $\left\lVert \cdot \right\rVert_1$ and $\left\lVert \cdot \right\rVert_\infty$ appearing in the objective function. A smooth version is easily obtained by rewritting both regularization terms using new decision variables. Since the first regularization term includes absolute values,
\begin{equation}
\nonumber
\left\lVert \bm{a}_{j\cdot}\right\rVert_1 = \sum_{t\in\tau_B} \left\lvert a_{jt}\right\rvert,\,\, j=1,\ldots,p,
\end{equation}
decision variables $a_{jt}\in\left[ -1,1\right],\,\, j=1,\ldots,p,\,\, t\in\tau_B$, are split into their positive and negative counterparts $a^+_{jt},\,\,a^-_{jt}\in\left[ 0,1\right],\,\, j=1,\ldots,p,\,\, t\in\tau_B$, respectively, holding $a_{jt} = a^+_{jt} - a^-_{jt}$ and $\left\lvert a_{jt}\right\rvert = a^+_{jt} + a^-_{jt}$. Similarly, we denote $\bm{a}^+ = \left( a^+_{jt}\right)_{j=1,\ldots,p,\,\,t\in\tau_B}$ and $\bm{a}^- = \left( a^-_{jt}\right)_{j=1,\ldots,p,\,\,t\in\tau_B}$. Regarding the second regularization term, new decision variables $\beta_j\in\left[0,1\right]$ are needed:
\begin{equation}
\nonumber
\left\lVert \bm{a}_{j\cdot} \right\rVert_\infty = \max_{t\in\tau_B} \left\lvert a_{jt}\right\rvert = \beta_j\in \left[0,1\right],\,\, j = 1,\ldots,p,
\end{equation}
and have to force $\beta_j \geq \left\lvert a_{jt}\right\rvert = a^+_{jt}+a^-_{jt},\,\, j =1,\ldots,p,\,\, t\in\tau_B$.

We can now formulate S-ORCT as a smooth problem, thus solvable
with standard continuous optimization solvers, as done in our computational section. Indeed, we have that \eqref{objetivo}-\eqref{res6} is equivalent to
\begingroup
\allowdisplaybreaks
\begin{align}
\min \hspace*{0.6cm}&  g\left(\bm{a}^+ - \bm{a}^-,\bm{\mu},\bm{C}\right) +\lambda^L \sum_{j=1}^p \sum_{t\in\tau_B} \left( a_{jt}^+ + a_{jt}^-\right) +\lambda^G \sum_{j=1}^p \beta_j \label{objetivocontinuo}\\
\text{s.t. \hspace*{0.5cm}}
& \sum_{k=1}^K C_{kt} = 1, \,\, t\in \tau_L,\label{res2c}\\
     & \sum_{t\in \tau_L} C_{kt} \geq 1, \,\, k=1,\ldots, K,\label{res3c} \\
     & \beta_j \geq a_{jt}^+ + a_{jt}^-,\,\, j=1,\ldots,p,\label{resnew}\\
     & a_{jt}^+,\,\, a_{jt}^- \in \left[ 0,1\right],\,\, j=1,\ldots,p,\,\, t\in \tau_B,\label{res4c}\\
     & \beta_j \in \left[0,1\right],\,\, j=1,\ldots,p,\label{resnew3}\\
     & \mu_{t} \in \left[-1,1\right] , \,\, t\in \tau_B,\label{res5c}\\
     & C_{kt} \in \left[ 0,1\right], \,\, k=1,\ldots,K,\,\, t\in \tau_L.\label{res6c}
\end{align}
\endgroup

Observe that, if we are only concerned about global sparsity, and thus we set $\lambda^L=0$, the rewriting of the decision variables $a_{jt},\,\, j=1,\ldots,p,\,\,t\in\tau_B$ is no longer necessary and \eqref{res4} replaces \eqref{res4c}, and \eqref{resnew} turns into  \begin{align}
\label{beta1}
\beta_j &\geq a_{jt}, \,\, j=1,\ldots,p,\,\,t\in\tau_B,\\
\label{beta2}
\beta_j &\geq -a_{jt}, \,\, j=1,\ldots,p,\,\,t\in\tau_B.
\end{align}


\section{Theoretical properties}
\label{Theoretical properties}

This section discusses some theoretical properties enjoyed by the S-ORCT. Let us consider the objective function of \eqref{objetivo}-\eqref{res6}. When taking $\lambda^L$ and $\lambda^G$ large enough, the first term related to the performance of the classifier becomes negligible and therefore $\bm{a}$ will shrink to $\bm{0}$. The tree with $\bm{a}=\bm{0}$ is the sparsest possible tree though not the best promising one from the accuracy point of view, since none of the predictor variables is used to classify. In this case, the probability of an individual with predictor variables $\bm{x}$ being assigned to class $k$ is independent of $\bm{x}$, and nothing more than the distribution of classes is available.
In this section, we derive upper bounds for the sparsity parameters, $\lambda^L$ and $\lambda^G$, in the sense that above these bounds the sparsest tree (with $\bm{a}^*=\bm{0}$) is a stationary point of the S-ORCT, that is, there exists  $\left(\bm{a}^*=\bm{0},\bm{\mu}^*,\bm{C}^*\right)$ such that the necessary optimality condition with respect to $\bm{a}$ is satisfied. This is done in Theorems \ref{globarsparsitytheorem} and \ref{teorema1}.

\begin{theo}
\label{globarsparsitytheorem}
Let $\sigma\in\left[0,1\right]$. For
\begin{align*}
\lambda^L \geq & \left(1-\sigma\right)\max_{\substack{\bm{\mu}\in\left[-1,1\right]^{\left\lvert\tau_B\right\rvert}\\ \bm{C}\in\left\lbrace 0,1\right\rbrace^{K\times\left\lvert\tau_L\right\rvert} }} \max_{j=1,\ldots,p} \left\lVert \nabla_{\bm{a}_{j\cdot}} g\left( \bm{0},\bm{\mu},\bm{C} \right)\right\rVert_\infty \text{and}\\
\lambda^G \geq & \hspace{0.5cm}\sigma \hspace{0.5cm}\max_{\substack{\bm{\mu}\in\left[-1,1\right]^{\left\lvert\tau_B\right\rvert}\\ \bm{C}\in\left\lbrace 0,1\right\rbrace^{K\times\left\lvert\tau_L\right\rvert} }} \max_{j=1,\ldots,p} \left\lVert \nabla_{\bm{a}_{j\cdot}} g\left( \bm{0},\bm{\mu},\bm{C} \right)\right\rVert_1,
\end{align*} $\bm{a}^*=\bm{0}$ is a stationary point of the S-ORCT.
\end{theo}

\noindent\textit{Proof.}

Let $\sigma,\,\,\lambda^L,\,\,\lambda^G$ be such that they satisfy the assumptions.

By Theorem \ref{resultadosoluciondeterminista}, there exists $\left(\bm{a}^*,\bm{\mu}^*,\bm{C}^* \right)$ optimal solution to \eqref{objetivo}-\eqref{res6} satisfying $C^*_{kt}\in\left\lbrace 0,1\right\rbrace \forall k=1,\ldots,K,\,\, t\in \tau_L$. In the following we will show that $\left(\bm{0},\bm{\mu}^*,\bm{C}^*\right)$ is a stationary point of the S-ORCT, i.e.,
\begin{equation}
\label{condicionnecesariaoptimalidad}
-\nabla_{\bm{a}} g\left( \bm{0},\bm{\mu}^*,\bm{C}^*\right)\in \partial_{\bm{a}}\left(\lambda^L\sum_{j=1}^p \left\lVert \bm{a}_{j\cdot} \right\rVert_1 + \lambda^G\sum_{j=1}^p \left\lVert \bm{a}_{j\cdot} \right\rVert_\infty\right)\left( \bm{0}\right)
\end{equation} where $\partial_{\bm{a}}$ is the subdifferential operator.

For every $\bm{a}_{j\cdot},\,\, j=1,\ldots,p$, we have that
\begin{align*}
\partial_{\bm{a}_{j\cdot}}\left(\left\lVert \bm{a}_{j\cdot} \right\rVert_1\right)\left(\bm{0}\right) = \mathbb{B}_\infty=\left\lbrace \bm{q}\in\mathbb{R}^{\left\lvert \tau_B\right\rvert}: \left\lVert \bm{q}\right\rVert_\infty \leq 1 \right\rbrace\\
\partial_{\bm{a}_{j\cdot}}\left(\left\lVert \bm{a}_{j\cdot} \right\rVert_\infty\right)\left(\bm{0}\right) = \mathbb{B}_1=\left\lbrace \bm{q}\in\mathbb{R}^{\left\lvert \tau_B\right\rvert}: \left\lVert \bm{q}\right\rVert_1 \leq 1 \right\rbrace.
\end{align*} Hence,
\begin{equation}
\nonumber
-\nabla_{\bm{a}_{j\cdot}} g\left( \bm{0},\bm{\mu}^*,\bm{C}^*\right)\in \lambda^L \partial_{\bm{a}_{j\cdot}}\left(\left\lVert \bm{a}_{j\cdot} \right\rVert_1\right)\left(\bm{0}\right) + \lambda^G \partial_{\bm{a}_{j\cdot}}\left(\left\lVert \bm{a}_{j\cdot} \right\rVert_\infty\right)\left(\bm{0}\right),
\end{equation}
 if, and only if,
 \begin{equation}
 \nonumber
 -\nabla_{\bm{a}_{j\cdot}} g\left( \bm{0},\bm{\mu}^*,\bm{C}^*\right)\in \lambda^L \mathbb{B}_\infty + \lambda^G \mathbb{B}_1,
\end{equation} if, and only if, there exist $\bm{q}^L_j,\,\, \bm{q}^G_j\in\mathbb{R}^{\left\lvert\tau_B\right\rvert}$ such that
\begin{align*}
&\left\lVert \bm{q}_j^L\right\rVert_\infty \leq 1,\\
&\left\lVert \bm{q}_j^G\right\rVert_1 \leq 1,\\
&-\nabla_{\bm{a}_{j\cdot}} g\left( \bm{0},\bm{\mu}^*,\bm{C}^*\right) = \lambda^L \bm{q}^L_j + \lambda^G \bm{q}^G_j,
\end{align*}
if, and only if, there exist $\tilde{\bm{q}}^L_j,\,\, \tilde{\bm{q}}^G_j\in\mathbb{R}^{\left\lvert\tau_B\right\rvert}$ such that
\begin{align*}
& \left\lVert \tilde{\bm{q}}_j^L\right\rVert_\infty \leq \lambda^L,\\
& \left\lVert \tilde{\bm{q}}_j^G\right\rVert_1 \leq \lambda^G,\\
& -\nabla_{\bm{a}_{j\cdot}} g\left( \bm{0},\bm{\mu}^*,\bm{C}^*\right) =  \tilde{\bm{q}}^L_j +  \tilde{\bm{q}}^G_j.
\end{align*} Let us consider
\begin{align*}
\tilde{\bm{q}}^L_j &=- \left( 1-\sigma\right) \nabla_{\bm{a}_{j\cdot}} g\left( \bm{0},\bm{\mu}^*,\bm{C}^*\right),\\
\tilde{\bm{q}}^G_j &=- \hspace{0.5cm}\sigma \hspace{0.5cm} \nabla_{\bm{a}_{j\cdot}} g\left( \bm{0},\bm{\mu}^*,\bm{C}^*\right),
\end{align*} and check that the conditions are satisfied:
\begin{align*}
& \left\lVert \tilde{\bm{q}}_j^L\right\rVert_\infty = \left( 1-\sigma\right) \left\lVert\nabla_{\bm{a}_{j\cdot}} g\left( \bm{0},\bm{\mu}^*,\bm{C}^*\right)\right\rVert_\infty \leq \left(1-\sigma\right)\max_{\substack{\bm{\mu}\in\left[-1,1\right]^{\left\lvert\tau_B\right\rvert}\\ \bm{C}\in\left\lbrace 0,1\right\rbrace^{K\times\left\lvert\tau_L\right\rvert} }} \max_{j=1,\ldots,p} \left\lVert \nabla_{\bm{a}_{j\cdot}} g\left( \bm{0},\bm{\mu},\bm{C} \right)\right\rVert_\infty \leq \lambda^L,\\
& \left\lVert \tilde{\bm{q}}_j^G\right\rVert_1 = \sigma \left\lVert\nabla_{\bm{a}_{j\cdot}} g\left( \bm{0},\bm{\mu}^*,\bm{C}^*\right)\right\rVert_1 \leq \sigma\max_{\substack{\bm{\mu}\in\left[-1,1\right]^{\left\lvert\tau_B\right\rvert}\\ \bm{C}\in\left\lbrace 0,1\right\rbrace^{K\times\left\lvert\tau_L\right\rvert} }} \max_{j=1,\ldots,p} \left\lVert \nabla_{\bm{a}_{j\cdot}} g\left( \bm{0},\bm{\mu},\bm{C} \right)\right\rVert_1\leq\lambda^G,\\
& \tilde{\bm{q}}^L_j +  \tilde{\bm{q}}^G_j = - \left( 1-\sigma\right) \nabla_{\bm{a}_{j\cdot}} g\left( \bm{0},\bm{\mu}^*,\bm{C}^*\right) - \sigma  \nabla_{\bm{a}_{j\cdot}} g\left( \bm{0},\bm{\mu}^*,\bm{C}^*\right) = -\nabla_{\bm{a}_{j\cdot}} g\left( \bm{0},\bm{\mu}^*,\bm{C}^*\right).
\end{align*}
Therefore, the desired result follows.  \qed

 A stronger result is proven for the S-ORCT of depth $D=1$ and $K=2$. Since local and global sparsity are equivalent for the S-ORCT of depth $D=1$, without loss of generality, we can assume that $\lambda^G=0$. Therefore, the objective function of the S-ORCT of depth $D=1$ can be written as:
\begin{equation}
\nonumber
g_1\left(\bm{a}_{\cdot 1},\mu_1,\bm{C}\right)=g\left(\bm{a}_{\cdot 1},\mu_1,\bm{C}\right) + \lambda^L \left\lVert \bm{a}_{\cdot 1}\right\rVert_1,
\end{equation}

where

\begin{align}
g\left(\bm{a}_{\cdot 1},\mu_1,\bm{C}\right) & = \dfrac{1}{N}\sum_{i=1}^N\left[ p_{i1}\left( \bm{a}_{\cdot 1},\mu_1\right) \sum_{k=1}^2 W_{y_ik}C_{k2} + \left( 1-p_{i1}\left( \bm{a}_{\cdot 1},\mu_1\right)\right) \sum_{k=1}^2 W_{y_i k}C_{k3} \right]\nonumber\\
  & = \dfrac{1}{N} \sum_{k=1}^2 \sum_{i\in I_k} \left[ p_{i1}\left( \bm{a}_{\cdot 1},\mu_1\right)\sum_{k'\neq k} W_{kk'}C_{k'2} + \left( 1-p_{i1}\left( \bm{a}_{\cdot 1},\mu_1\right)\right) \sum_{k'\neq k} W_{kk'}C_{k'3}\right] \label{gdepth1}
\end{align} and
\begin{equation}
\nonumber
p_{i1}\left( \bm{a}_{\cdot 1},\mu_1\right)=F\left(\dfrac{1}{p}\bm{a}_{\cdot 1}^T \bm{x}_i-\mu_1\right),\,\, i=1,\ldots,N.
\end{equation}

A technical lemma is needed to prove the desired result.

\begin{lemma}\label{lemma1} For any allocation rule $\bm{C}$, the objective function of the S-ORCT of depth $D=1$, $g_1$, is monotonic in $\mu_1$ when $\bm{a}_{\cdot 1}=\bm{0}$.  \end{lemma}

\noindent\textit{Proof.}

Fixed $\bm{a}_{\cdot 1}=\left( a_{j1}\right)_{j=1,\ldots,p}$, and $\bm{C}=\left( C_{kt}\right)_{k=1,2,\,\, t=2,3}$,
\begin{equation}
\nonumber
\left.\dfrac{\partial g_1}{\partial \mu_1}\right|_{\bm{a}_{\cdot 1}=\bm{0}} = \dfrac{1}{N}\sum_{k=1}^K \sum_{i\in I_k} \left(\sum_{k'\neq k}W_{kk'} C_{k'2} -\sum_{k'\neq k} W_{kk'}C_{k'3}\right) \left.\dfrac{\partial  p_{i1}\left( \bm{a}_{\cdot 1},\mu_1\right)}{\partial \mu_1}\right|_{\bm{a}_{\cdot 1}=\bm{0}},
\end{equation} where \begin{align*}
\nonumber
\dfrac{\partial  p_{i1}\left( \bm{a}_{\cdot 1},\mu_1\right)}{\partial \mu_1} &= \dfrac{\partial  F\left( \dfrac{1}{p}\bm{a}_{\cdot 1}^T \bm{x}_i-\mu_1\right)}{\partial  \left(\dfrac{1}{p}\bm{a}_{\cdot 1}^T \bm{x}_i-\mu_1\right)}\dfrac{\partial  \left(\dfrac{1}{p}\bm{a}_{\cdot 1}^T \bm{x}_i-\mu_1\right)}{\partial \mu_1} = -f \left(\dfrac{1}{p}\bm{a}_{\cdot 1}^T \bm{x}_i-\mu_1\right), \,\, i=1,\ldots,N,
\end{align*} and \begin{align*}
\nonumber
\left.\dfrac{\partial  p_{i1}\left( \bm{a}_{\cdot 1},\mu_1\right)}{\partial \mu_1}\right|_{\bm{a}_{\cdot 1}=\bm{0}} = -f\left( -\mu_1\right),\,\, i=1,\ldots,N.
\end{align*} Thus, \begin{align*}
\nonumber
\left.\dfrac{\partial g_1\left( \bm{a}_{\cdot 1},\mu_1,\bm{C}\right)}{\partial \mu_1}\right|_{\bm{a}_{\cdot 1}=\bm{0}} &= \dfrac{1}{N} f\left(-\mu_1\right)\left(\sum_{i\in I_1} W_{12}\left( C_{23}-C_{22}\right)+\sum_{i\in I_2} W_{21}\left( C_{13}-C_{12}\right)\right)\\
&= \dfrac{1}{N} f\left( -\mu_1\right)\left( W_{12}\left( C_{23}-C_{22}\right)\lvert I_1\rvert + W_{21}\left( 1-C_{23}-1+C_{22}\right)\lvert I_2\rvert\right)\\
&= \dfrac{1}{N} f\left( -\mu_1\right)\left( C_{23}-C_{22}\right)\left( W_{12}\lvert I_1\rvert - W_{21}\lvert I_2\rvert\right).
\end{align*}  Since $f$ is a probability density function, the expression $\left.\dfrac{\partial g_1\left( \bm{a}_{\cdot 1},\mu_1,\bm{C}\right)}{\partial \mu_1}\right|_{\bm{a}_{\cdot 1}=\bm{0}}$ will always have the same sign for any value of $\mu_1$ and the desired result follows. \qed

\begin{theo}\label{teorema1}
For 

\begin{align}
\label{lambdateorema3}
\lambda^L &\geq \dfrac{1}{N} \max_{j=1,\ldots,p} \left| -W_{21}\sum_{i\in I_2}x_{ij} + W_{12} \sum_{i\in I_1}x_{ij}\right| \max_{\mu_1\in\left\lbrace -1,1\right\rbrace} f\left( \mu_1\right),
\end{align}
 $\bm{a}_{\cdot 1}^*=\bm{0}$ is a stationary point of the S-ORCT of depth $D=1$.

\end{theo}

\noindent\textit{Proof.}

Using the monotonicity of $\mu_1$ proven in Lemma \ref{lemma1} and Theorem \ref{globarsparsitytheorem} with $\sigma =0$, we have that for
\begin{align}
\nonumber
\lambda^L &\geq \max_{\substack{\mu_1\in\left\lbrace -1,1\right\rbrace\\ \bm{C}\in\left\lbrace 0,1\right\rbrace^{2\times2} }} \max_{j=1,\ldots,p} \left\lvert \nabla_{\bm{a}_{j1}} g\left( \bm{0},\mu_1,\bm{C} \right)\right\rvert\\
\label{lambdausandoteorema1}
&= \max_{\substack{\mu_1\in\left\lbrace -1,1\right\rbrace\\ \bm{C}\in\left\lbrace 0,1\right\rbrace^{2\times2} }} \left\lVert \nabla_{\bm{a}_{\cdot 1}} g\left( \bm{0},\mu_1,\bm{C} \right)\right\rVert_\infty,
\end{align} where $g$ is as in \eqref{gdepth1}, $\bm{a}_{\cdot 1}^*=\bm{0}$ is a stationary point of thr S-ORCT. The remainder of the proof is devoted to rewriting \eqref{lambdausandoteorema1} as in \eqref{lambdateorema3}.

We proceed with the calculation of the gradient.

\noindent For $j =1,\ldots,p$:
\begin{equation}
\nonumber
\dfrac{\partial g \left(\bm{0},\mu_1,\bm{C}\right)}{\partial a_{j1}} = \left.\dfrac{\partial g \left(\bm{a}_{\cdot 1},\mu_1,\bm{C}\right)}{\partial a_{j1}}\right|_{\bm{a}_{\cdot 1}=\bm{0}} = \dfrac{1}{N} \sum_{k=1}^2 \sum_{i\in I_k} \left(\sum_{k'\neq k} W_{kk'}C_{k'2} -\sum_{k'\neq k} W_{kk'}C_{k'3}\right) \left.\dfrac{\partial  p_{i1}\left(\bm{a}_{\cdot 1},\mu_1 \right)}{\partial a_{j1}}\right|_{\bm{a}_{\cdot 1}=\bm{0}},
\end{equation} where
\begin{align*}
\dfrac{\partial  p_{i1}\left( \bm{a}_{\cdot 1},\mu_1\right)}{\partial a_{j1}} &= \dfrac{\partial  F\left( \dfrac{1}{p}\bm{a}_{1}^T \bm{x}_i-\mu_1\right) }{\partial  \left(\dfrac{1}{p}\bm{a}_{1}^T \bm{x}_i-\mu_1\right)}\dfrac{\partial \left( \dfrac{1}{p}\bm{a}_{1}^T \bm{x}_i-\mu_1\right)}{\partial a_{j1}}= \dfrac{x_{ij}}{p}f\left( \dfrac{1}{p}\bm{a}_{1}^T \bm{x}_i-\mu_1\right) , \,\, i=1,\ldots,N.
\end{align*} and
\begin{equation}
\nonumber
\left.\dfrac{\partial  p_{i1}\left( \bm{a}_{\cdot 1},\mu_1\right)}{\partial a_{j1}}\right|_{\bm{a}_{\cdot 1}=\bm{0}} = \dfrac{x_{ij}}{p}f\left( -\mu_1\right) , \,\, i=1,\ldots,N.
\end{equation} Thus,
\begin{align*}
\dfrac{\partial g \left(\bm{0},\mu_1,\bm{C}\right)}{\partial a_{j1}} & = \dfrac{1}{Np}f\left( -\mu_1\right)\left( W_{12}\sum_{i\in I_1} x_{ij}\left( C_{22}-C_{23}\right) + W_{21}\sum_{i\in I_2} x_{ij}\left( C_{12}-C_{13}\right) \right).\label{derivada parcial}
\end{align*}

Now, we look for the maximum $\lambda^L$ among every possible allocation of the decision variables $\bm{C}$, i.e.:
\begin{equation}
\nonumber
\lambda^L_{\mu_1}=\max_{\bm{C}\in\left\lbrace 0,1\right\rbrace^{2\times2}}\left\lVert \nabla_{\bm{a}_{\cdot 1}} g\left( \bm{0},\mu_1,\bm{C} \right)\right\rVert_\infty=\max_{\overline{\bm{C}}\in\left\lbrace 0,1\right\rbrace^{4\times 1}} \,\, \lVert D\overline{\bm{C}}\rVert_{\infty},
\end{equation}
where
\begin{equation}
\nonumber
D=\dfrac{1}{Np}f\left( -\mu_1\right)\begin{pmatrix}
-W_{21}\sum_{i\in I_2}x_{i1} & W_{21}\sum_{i\in I_2}x_{i1} & -W_{12}\sum_{i\in I_1}x_{i1}&W_{12}\sum_{i\in I_1}x_{i1}\\
\vdots & \vdots & \vdots & \vdots\\
-W_{21}\sum_{i\in I_2}x_{ip} & W_{21}\sum_{i\in I_2}x_{ip} & -W_{12}\sum_{i\in I_1}x_{ip}&W_{12}\sum_{i\in I_1}x_{ip}
\end{pmatrix}
\end{equation} and $\overline{\bm{C}}=\left( C_{12},C_{13}, C_{22},C_{23} \right)^T$.
\begin{align*}
\max_{\overline{\bm{C}}\in\left\lbrace 0,1\right\rbrace^{4\times 1}} \,\, \lVert D\overline{\bm{C}}\rVert_{\infty} &= \max_{\overline{\bm{C}}\in\left\lbrace 0,1\right\rbrace^{4\times 1}} \,\, \max \,\,\left\lbrace \vert d_1^{T}\,\,\overline{\bm{C}}\vert,\ldots, \vert d_p^{T}\,\,\overline{\bm{C}}\vert \right\rbrace\\
& = \max_{\overline{\bm{C}}\in\left\lbrace 0,1\right\rbrace^{4\times 1}} \,\, \max\,\,\left\lbrace d_1^{T}\,\,\overline{\bm{C}}, -d_1^{T}\,\,\overline{\bm{C}}, \ldots, d_p^{T}\,\,\overline{\bm{C}}, -d_p^{T}\,\,\overline{\bm{C}} \right\rbrace\\
& = \max \left\lbrace \max_{\overline{\bm{C}}\in\left\lbrace 0,1\right\rbrace^{4\times 1}}\,\, d_1^{T}\,\,\overline{\bm{C}}, \max_{\overline{\bm{C}}\in\left\lbrace 0,1\right\rbrace^{4\times 1}}\,\, -d_1^{T}\,\,\overline{\bm{C}},\ldots, \max_{\overline{\bm{C}}\in\left\lbrace 0,1\right\rbrace^{4\times 1}}\,\, d_p^{T}\,\,\overline{\bm{C}}, \max_{\overline{\bm{C}}\in\left\lbrace 0,1\right\rbrace^{4\times 1}}\,\, -d_p^{T}\,\,\overline{\bm{C}}  \right\rbrace.
\end{align*}

A finite number of transportation problems is to be solved, with the form:
\begin{align*}
z=\max_{\overline{\bm{C}}\in\left\lbrace 0,1\right\rbrace^{4\times 1}} \,\, & \left\lbrace \pm d_j^T\,\,\overline{\bm{C}}\right\rbrace\\
\,\,\text{s.t.} \,\,\,\,\,& C_{12}+C_{22} = 1\\
& C_{13}+C_{23} = 1\\
& C_{12}+C_{13} \geq 1\\
& C_{22}+C_{23} \geq 1,
\end{align*}
for which the integrality property holds. Then, we only have as possible solutions: $\overline{\bm{C}}=\left(1,0,0,1\right)^T$ or $\overline{\bm{C}}=\left(0,1,1,0\right)^T$. Thus, the optimal objective is obtained as follows:
\begin{align*}
z_{\textit{opt}}&=\max\left\lbrace \left. \pm d_j^T\,\, \overline{\bm{C}}\right|_{\overline{\bm{C}}=\left( 1,0,0,1\right)^T},\left. \pm d_j^T\,\, \overline{\bm{C}}\right|_{\overline{\bm{C}}=\left( 0,1,1,0\right)^T} \right\rbrace\\
& = \max\left\lbrace \dfrac{1}{Np}f\left( -\mu_1\right) \left(-W_{21}\sum_{i\in I_2}x_{ij} + W_{12}\sum_{i\in I_1}x_{ij}\right), \dfrac{1}{Np}f\left( -\mu_1\right)\left(W_{21}\sum_{i\in I_2}x_{ij} - W_{12}\sum_{i\in I_1}x_{ij}\right)\right\rbrace\\
& = \dfrac{1}{Np}f\left( -\mu_1\right) \left| -W_{21}\sum_{i\in I_2}x_{ij} + W_{12}\sum_{i\in I_1}x_{ij} \right|.
\end{align*}

Let us define
\begin{equation}
\nonumber
\lambda^L_{\mu_1} = \dfrac{1}{Np}f\left( -\mu_1\right) \max_{j=1,\ldots,p} \left| -W_{21}\sum_{i\in I_2}x_{ij} + W_{12} \sum_{i\in I_1}x_{ij}\right|,
\end{equation}
and the result holds when
\begin{align*}
\lambda^L &\geq \max\left\lbrace \lambda^L_{\mu_1=-1},\lambda^L_{\mu_1=1}\right\rbrace.
\end{align*}\qed

\section{Computational experience}
\label{Computational experience}

\subsection{Introduction}

The aim of this section is to illustrate the performance of our sparse optimal randomized classification trees S-ORCT's. We have run our model for a grid of values of the sparsity regularization parameters $\lambda^L$ and $\lambda^G$. The message that can be drawn from our experimental experience is twofold. First, we show empirically that our S-ORCT can gain in both local and global sparsity, without harming classification accuracy. Second, we benchmark our approach against CART, the classic approach to build decision trees, which considers orthogonal cuts and therefore has the best possible local sparsity. We show that we are able to trade in some of our classification accuracy, still being superior to CART, to be comparable to CART in terms of global sparsity.

The S-ORCT smooth formulation \eqref{objetivocontinuo}-\eqref{res6c} has been implemented using Pyomo optimization modeling language \cite{hart2017pyomo,hart2011pyomo} in Python 3.5 \cite{pthn}. As solver, we have used IPOPT 3.11.1 \cite{Wachter2006}, and have followed a multistart approach, where the process is repeated $20$ times starting from different random initial solutions. For CART, the implementation in the \texttt{rpart} R package \cite{rpartR} is used. Our experiments have been conducted on a PC, with an Intel$^\circledR$ Core$^{\rm TM}$ i7-2600 CPU 3.40GHz processor and 16 GB RAM. The operating system is 64 bits.

 The remainder of the section is structured as follows. Section \ref{setup} gives details on the procedure followed to test S-ORCT. In Sections \ref{Results for local S-ORCTs} and \ref{Results for global sparsity}, respectively, we discuss the results for local and global sparsities separately, while in Section \ref{Results for local and global sparsity} we present results when both sparsities are simultaneously taken into account. Finally, Section \ref{Comparison S-ORCT versus CART} statistically compares S-ORCT versus CART in terms of classification accuracy and global sparsity.

\subsection{Setup}
\label{setup}

An assorted collection of well-known real data sets from the UCI Machine Learning Repository \cite{Lichman:2013} has been chosen for the computational experiments. Table \ref{Information about the data sets considered.} lists their names together with their number of observations, number of predictor variables and number of classes with the corresponding class distribution. In our pursuit of building small and, therefore, less complex trees, the construction of S-ORCTs has been restricted to depth $D=1$ for two-class problems and depth $D=2$ for three- and four- class problems.

\begin{table}[h!] 
\caption{Information about the data sets considered.}
\label{Information about the data sets considered.}
\centering
\begin{tabular}{| l | l | c | c | c | c |}
  \hline			
  Data set & Abbrev. & $N$ & $p$ & $K$ & Class distribution\\ \hline
Monks-problems-3 & Monks-3 &122 & 11 & 2 & 51\% - 49\% \\
\hline
Monks-problems-1 & Monks-1 & 124 & 11 & 2 & 50\% - 50\% \\
\hline
Monks-problems-2 & Monks-2 & 169 &11 & 2 & 62\% - 38\% \\
\hline
Connectionist-bench-sonar & Sonar & 208 & 60 & 2 & 55\% - 45\% \\
  \hline
  Ionosphere & Ionosphere & 351 & 34 & 2 & 64\% - 36\%\\
  \hline
  Breast-cancer-Wisconsin & Wisconsin & 569 & 30 & 2 &63\% - 37\%\\
  \hline
  Credit-approval & Creditapproval& 653 & 37 & 2 & 55\% - 45\% \\
  \hline
  Pima-indians-diabetes & Pima &768 & 8& 2& 65\% - 35\% \\
  \hline
 Statlog-project-German-credit & Germancredit& 1000 & 48 & 2 & 70\% - 30\% \\
  \hline
  Banknote-authentification & Banknote & 1372 & 4 & 2 & 56\% - 44\% \\
  \hline
  Ozone-level-detection-one & Ozone & 1848 & 72 & 2 & 97\% - 3\% \\
  \hline
 Spambase & Spam & 4601 & 57 & 2& 61\% - 39\% \\
  \hline
  Iris & Iris & 150 & 4 & 3 & 33.3\%-33.3\%-33.3\% \\
  \hline
    Wine & Wine & 178 & 13 & 3 & 40\%-33\%-27\% \\
  \hline
    Seeds & Seeds & 210 & 7 & 3 & 33.3\%-33.3\%-33.3\% \\
  \hline
 Balance-scale & Balance & 625& 16 &3 & 46\%-46\%-8\% \\
 \hline
  Thyroid-disease-ann-thyroid & Thyroid & 3772 & 21 & 3 & 92.5\%-5\%-2.5\% \\
  \hline
  Car-evaluation & Car & 1728 & 15 & 4 & 70\%-22\%-4\%-4\% \\
  \hline
\end{tabular}
\end{table}
Each data set has been split into two subsets: the training subset (75\%) and the test subset (25\%). The corresponding S-ORCT is built on the training subset and, then, accuracy, local and global sparsities are measured. The out-of-sample accuracy over the test subset is denoted by \textit{acc}. Local sparsity is denoted by $\delta^L$ and reads as the average percentage of predictor variables not used per branch node:
\begin{align*}
\delta^L = \dfrac{1}{\left\lvert\tau_B\right\rvert}\sum_{t\in\tau_B}\dfrac{\left\lvert\left\lbrace a_{jt}= 0,\,\, j=1,\ldots,p\right\rbrace\right\rvert}{p}\times 100.
\end{align*} Global sparsity, $\delta^G$, is measured as the percentage of predictor variables not used at any of the branch nodes, i.e., across the whole tree:
\begin{align*}
\delta^G = \dfrac{\left\lvert \left\lbrace \bm{a}_{j\cdot}= \bm{0},\,\,j=1,\ldots,p\right\rbrace\right\rvert}{p}\times 100.
\end{align*}
Note that when $D=1$, local and global sparsity are measuring the same since there is a single cut across the whole tree. The training/testing procedure has been repeated ten times in order to avoid the effect of the initial split of the data. The results shown in the tables represent the average of such ten runs to each of the three performance criteria.

In what follows, we describe the choices made for the parameters in S-ORCT. Equal misclassification weights, $W_{y_ik} = 0.5,\,\, k=1,\ldots,K,\,\,k\neq y_i$, have been used for the experiments. We have added the set of constraints \eqref{performanceconstraints} with $\rho_k=0.1,\,\,k=1,\ldots,K$. The logistic CDF has been chosen for our experiments: 
\begin{equation}
\nonumber
F\left(\cdot\right) = \dfrac{1}{1+\exp\left(-\left(\cdot\right)\gamma\right)},
\end{equation} 
with a large value of $\gamma$, namely, $\gamma = 512$. The larger the value of $\gamma$, the closer the decision rule defined by $F$ is to a deterministic rule. We will illustrate that a small level of randomization is enough for obtaining good results. We have trained S-ORCT, as formulated in \eqref{objetivocontinuo}-\eqref{res6c}, for $17 \times 17$ pairs of values for $\left(\lambda^L,\lambda^G\right)$ starting from $\lambda^L=0$ followed by the grid $\left\lbrace \dfrac{2^r}{p\left\lvert\tau_B\right\rvert},\,\, -12\leq r\leq 3,\,\, r\in\mathbb{Z} \right\rbrace$, and, similarly, $\lambda^G=0$ followed by the grid $\left\lbrace \dfrac{2^r}{p},\,\, -12\leq r\leq 3,\,\, r\in\mathbb{Z} \right\rbrace$. We start solving the optimization problem with $\left(\lambda^L,\lambda^G\right)=\left(0,0\right)$, where the multistart approach uses 20 random initial solutions. We continue solving the optimization problem for $\lambda^L=0$ but with larger values of $\lambda^G$. Once all values of $\lambda^G$ are executed, we start the process all over again with the next value of $\lambda^L$ in the grid. For pair $\left(\lambda^L,\lambda^G\right)$, we feed the corresponding optimization problem with the 20 solutions resulting from the problem solved for the previous pair. For a given initial solution, the computing time taken by the S-ORCT typically ranges from 0.33 seconds (in Monks-1) to 22.27 seconds (in Thyroid).

For CART, the default parameter setting in \texttt{rpart} is used.

\subsection{Results for local sparsity}
\label{Results for local S-ORCTs}

\begin{figure} \hspace*{-1cm}
 \subfloat{
	\label{f:gato}
	\includegraphics[width=0.36\textwidth]{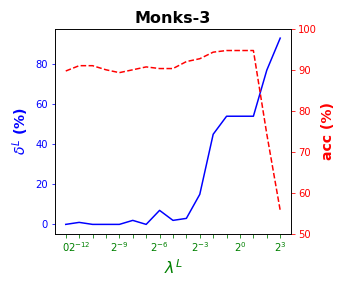}}
\subfloat{
	\label{f:tigre}
	\includegraphics[width=0.36\textwidth]{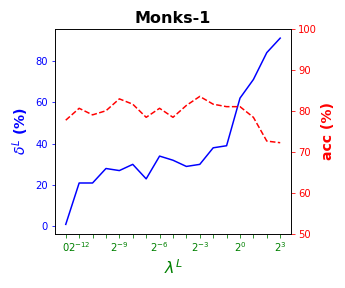}}
\subfloat{
	\label{f:conejo}
	\includegraphics[width=0.36\textwidth]{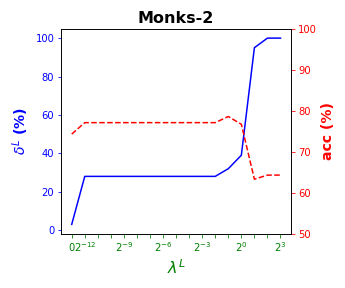}}\\
\hspace*{-1cm}  \subfloat{
   \label{f:gato}
    \includegraphics[width=0.36\textwidth]{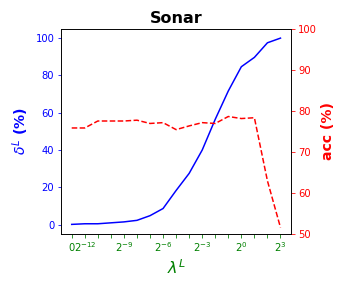}}
  \subfloat{
   \label{f:tigre}
    \includegraphics[width=0.36\textwidth]{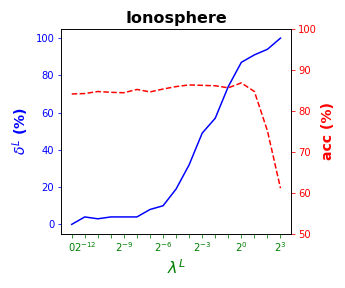}}
  \subfloat{
   \label{f:conejo}
    \includegraphics[width=0.36\textwidth]{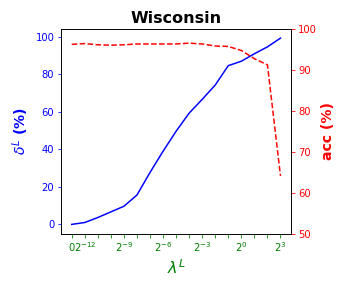}}\\
\hspace*{-1cm}\subfloat{
   \label{f:conejo}
    \includegraphics[width=0.36\textwidth]{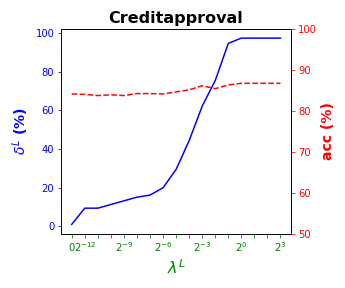}}
        \subfloat{
   \label{f:conejo}
    \includegraphics[width=0.36\textwidth]{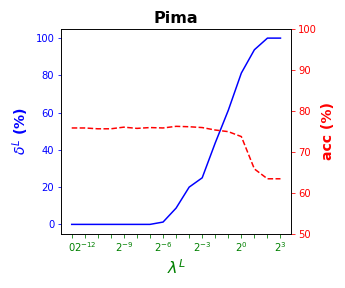}}
          \subfloat{
   \label{f:conejo}
    \includegraphics[width=0.36\textwidth]{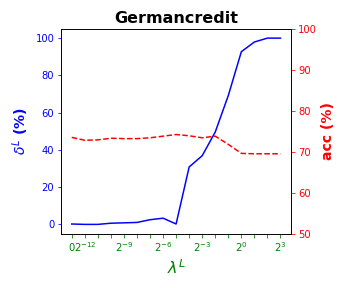}}\\
\hspace*{-1cm} \subfloat{
   \label{f:conejo}
    \includegraphics[width=0.36\textwidth]{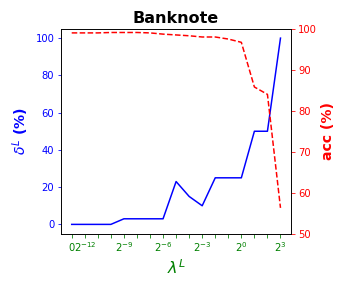}}
      \subfloat{
   \label{f:conejo}
    \includegraphics[width=0.36\textwidth]{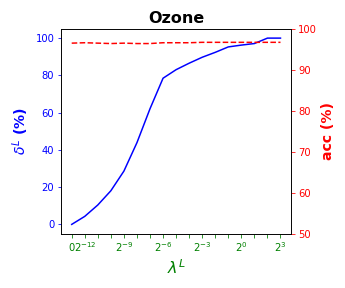}}
        \subfloat{
   \label{f:conejo}
    \includegraphics[width=0.36\textwidth]{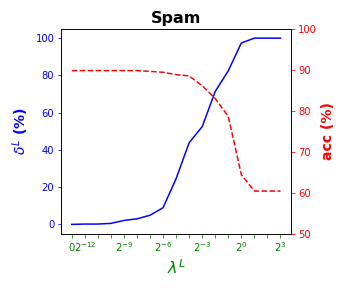}}
 \caption{ Graphical representation, for each data set, of the average percentage of predictor variables per branch node, $\delta^L$, together with the average out-of-sample accuracy obtained, \textit{acc}, as a function of the values of $\lambda^L$ considered in the local S-ORCT construction.}
 \label{f:animales}
\end{figure}

\begin{figure} \hspace*{-1cm}
	\subfloat{
		\label{f:gato}
		\includegraphics[width=0.36\textwidth]{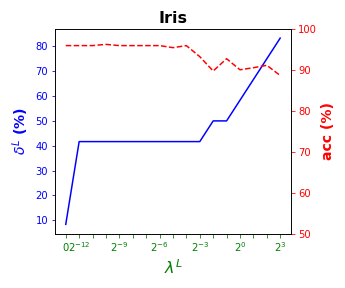}}
	\subfloat{
		\label{f:tigre}
		\includegraphics[width=0.36\textwidth]{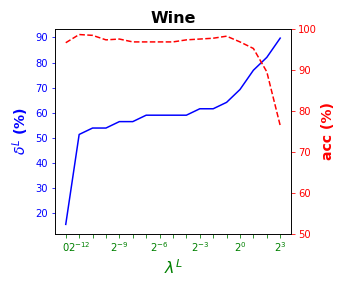}}
	\subfloat{
		\label{f:conejo}
		\includegraphics[width=0.36\textwidth]{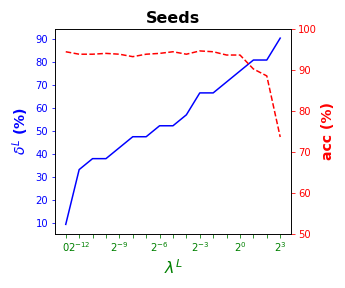}}\\
	\hspace*{-1cm}  \subfloat{
		\label{f:gato}
		\includegraphics[width=0.36\textwidth]{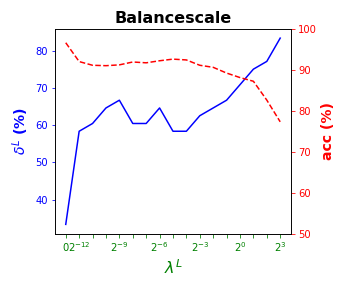}}
	\subfloat{
		\label{f:tigre}
		\includegraphics[width=0.36\textwidth]{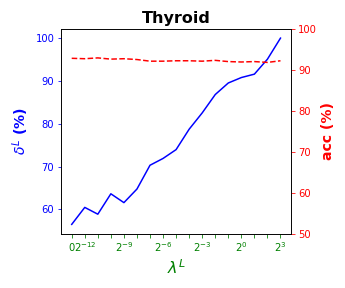}}
	\subfloat{
		\label{f:conejo}
		\includegraphics[width=0.36\textwidth]{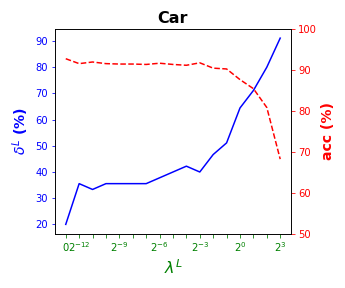}}
	\caption{ Graphical representation, for each data set, of the average percentage of predictor variables per branch node, $\delta^L$, together with the average out-of-sample accuracy obtained, \textit{acc}, as a function of the values of $\lambda^L$ considered in the local S-ORCT construction.}
	\label{f:animales22}
\end{figure}

Tables \ref{Results for $D=1$.} and \ref{Results} present the results of the so-called local S-ORCT, i.e., when $\lambda^G=0$ and thus only local sparsity is taken into account. Figures \ref{f:animales} and \ref{f:animales22} depict these results per data set, by showing simultaneously $\delta^L$ (blue solid line) and \textit{acc} (red dashed line) as a function of the grid of the $\lambda^L$'s considered. As expected, the larger the $\lambda^L$, the larger the $\delta^L$. The sparsest tree is shown in most of the data sets for large values of the parameter $\lambda^L$, where the best solution in terms of sparsity is obtained but the worst possible one in terms of accuracy. In terms of accuracy, the best rates are sometimes achieved when not all the predictor variables are included in the model. For instance, best performance is reached when sparsity is about $9-25\%$ for Pima, the $30\%$ for Monks-1, the $32\%$ for Monks-2, the $44\%$ for Germancredit, the $47\%$ for Car, the $52-56\%$ for Thyroid, the $54\%$ for Monks-3, the $55-60\%$ for Iris, the $72-90\%$ for Sonar, the $81\%$ for both Wine and Seeds and the $87\%$ for Ionosphere. We highlight the Creditapproval data set, on which one single predictor variable can already guarantee very good accuracy. For Ozone, accuracy remains over the $96\%$ for the grid of $\lambda^L$'s considered. Accuracy might be slightly damaged but a great gain in sparsity is obtained. This is the case for Banknote, Spam, Balance or Wisconsin, which present a loss of accuracy lower than the $1$ percentage point (p.p.), $4$ p.p., $6$ p.p. and $1$ p.p. but $25\%$, $52\%$, $63\%$ and $85\%$ of local sparsity is reached, respectively.
\begin{landscape}
\begin{table}[h!]
\caption{Results for the local S-ORCT of depth $D=1$ as a function of $\lambda^L$, where $\delta^L$ represents the average percentage of predictor variables not used per branch node in the tree over the ten runs and \textit{acc}, the average out-of-sample accuracy.}
\label{Results for $D=1$.}
\centering\hspace*{-0.5cm} {\small
\begin{tabular}{|c|c|c|c|c|c|c|c|c|c|c|c|c|c|c|c|c|c|c|c|c|c|c|c|c|}\hline
 \multirow{2}{*}{$\lambda^L$}  & \multicolumn{2}{|c|}{Monks-3} & \multicolumn{2}{|c|}{Monks-1}&\multicolumn{2}{|c|}{Monks-2} &  \multicolumn{2}{|c|}{Sonar} &  \multicolumn{2}{|c|}{Ionosphere} &\multicolumn{2}{|c|}{Wisconsin} &\multicolumn{2}{|c|}{Creditapproval}&\multicolumn{2}{|c|}{Pima}&\multicolumn{2}{|c|}{Germancredit}&  \multicolumn{2}{|c|}{Banknote}&\multicolumn{2}{|c|}{Ozone}&\multicolumn{2}{|c|}{Spam} \\\cline{2-25}
  & $\delta^L$ & acc & $\delta^L$ & acc & $\delta^L$ & acc & $\delta^L$ & acc &$\delta^L$ & acc & $\delta^L$ & acc & $\delta^L$ & acc & $\delta^L$ & acc & $\delta^L$ & acc & $\delta^L$ & acc & $\delta^L$ & acc & $\delta^L$ & acc \\\hline

$0$ &0 &89.7 &1 &77.7 & 3& 74.3 & 0 & 75.8 &0 &84.1 & 0 & 96.2 & 1 & 84.1& 0& 75.8& 0 & 73.5& 0 & 99.0& 0 & 96.5& 0 & 89.8\\ \hline
$2^{-12}$ &1 &91.0 &21 &80.6 &28 &77.1  & 0 & 75.8&4 & 84.2& 1 & 96.4 & 9&84.0 &0 &75.8 &0 &72.8 &0 & 99.0&4 & 96.6& 0&89.8 \\ \hline
$2^{-11}$ &0 &91.0 &21 &79.0 &28 &77.1  & 0 & 77.5&3 &84.7 & 4 & 96.1 & 9& 83.7& 0&75.6 & 0& 72.9&0 &99.0 & 10& 96.5&0 &89.8 \\ \hline
$2^{-10}$ & 0&90.0 &28 &80.0 &28 &77.1  & 1 & 77.5&4 &84.5 & 7 & 96.0 & 11&83.9 &0 &75.6 &1 &73.3 &0 &99.1 &18 &96.4 & 1&89.8 \\ \hline
$2^{-9}$ & 0&89.3 &27 &82.9 &28 &77.1  & 2 & 77.5& 4&84.4 & 10 & 96.1& 13 & 83.7& 0& 76&1 & 73.2& 3&99.1 & 29&96.5 &2 & 89.8\\ \hline
$2^{-8}$ &2 &90.0 &30 &81.6 &28 & 77.1 & 2 &77.7& 4&85.2 & 16&96.3 &15 &84.2 &0 &75.7 & 1&73.2 &3 &99.1 &44 &96.4 & 3&89.8 \\ \hline
$2^{-7}$ &0 &90.7 &23 &78.4 &28 &77.1  & 5&76.9 &8 &84.6 &28 &96.3 & 16& 84.2&0 &75.9 &2 & 73.4&3 &99.0 & 62&96.4 &5 & 89.6\\ \hline
$2^{-6}$ &7 &90.3 &34 &80.6 &28 &77.1  & 9 & 77.1 & 10&85.3 &39 & 96.3& 20 &84.1 &1 &75.8 &3 &73.8 &3 &98.7 & 78& 96.6& 9 &89.4 \\ \hline
$2^{-5}$ &2 &90.3 &32 &78.4 & 28& 77.1 & 18 & 75.4 & 19&85.9 &50 &96.3 & 29&84.6 & 9& 76.2& 0&74.2 &23 &98.5 & 83& 96.6&25 & 88.8\\ \hline
$2^{-4}$ &3 &92.0 &29 &81.3 &28 & 77.1 & 28 &76.3 &32 &86.3 & 59& 96.5&44 &85.1 &20 &76.1 & 31& 73.9& 15& 98.3& 87&96.6 & 44& 88.5\\ \hline
$2^{-3}$ &15 &92.7 &30 &83.5 & 28&77.1  & 40 & 77.1& 49&86.2 & 67& 96.3& 62&86.1 & 25&75.9 &37 & 73.4&10 & 98.0& 90&96.7 & 52& 86.1\\ \hline
$2^{-2}$ & 45&94.3 & 38&81.6 &28 &77.1 & 56 & 76.9& 57&86.1 & 74 & 95.8&75 &85.4 & 44& 75.3& 50&73.8 & 25&98.0 &92 &96.7 &71 & 83\\ \hline
$2^{-1}$ &54 &94.7 &39 &81.0 & 32&78.6  & 72 &78.6 & 74& 85.6& 85& 95.7& 95&86.3 & 61&74.9 &69 &71.8 &25 &97.5 &95 &96.7 & 82& 78.6\\ \hline
$2^{0}$ & 54&94.7 &62 &81.0 & 39& 76.7 & 85 &78.1& 87&86.8 &87 &94.7 &97 & 86.7& 81&73.7 &93 &69.6 & 25& 96.7& 96& 96.7&97 & 64.4\\ \hline
$2^{1}$ &54 &94.7 &71 &78.4 &95 & 63.3& 90 &78.3 &91 &84.7 & 91&92.7 & 97&86.7 &94 &65.8 &98 &69.5 &50 &85.8 &97 &96.7 & 100& 60.4\\ \hline
$2^{2}$ & 77&74.3 &84 &72.6 &100 &64.3  & 98 & 62.9 & 94&75.1 &95 &91.2 & 97&86.7 & 100& 63.4& 100& 69.5& 50& 84.0& 100& 96.7& 100& 60.4\\ \hline
$2^{3}$ &93 &55.7 &91 &72.2 &100 & 64.3 & 100 & 51.5& 100&61.1 & 99 & 64.1& 97& 86.7& 100& 63.4& 100& 69.5& 100& 56.3& 100& 96.7& 100& 60.4\\ \hline
\end{tabular}}
\end{table}
\end{landscape}

\begin{table}[h!]
	\caption{Results for the local S-ORCT of depth $D=2$ as a function of $\lambda^L$, where $\delta^L$ represents the average percentage of predictor variables not used per branch node in the tree over the ten runs and \textit{acc}, the average out-of-sample accuracy.}
	\label{Results}
	\centering
	\begin{tabular}{|c|c|c|c|c|c|c|c|c|c|c|c|c|}\hline
		\multirow{2}{*}{$\lambda^L$}  & \multicolumn{2}{|c|}{Iris} &\multicolumn{2}{|c|}{Wine} &\multicolumn{2}{|c|}{Seeds}&\multicolumn{2}{|c|}{Balance}&\multicolumn{2}{|c|}{Thyroid}&\multicolumn{2}{|c|}{Car} \\\cline{2-13}
		& $\delta^L$ & acc & $\delta^L$ & acc & $\delta^L$ & acc & $\delta^L$ & acc & $\delta^L$ & acc & $\delta^L$ & acc \\\hline
		$0$  & 8  & 95.9 & 15 & 96.6 &10 & 94.4 & 33 & 96.6& 57& 92.8& 20& 92.7 \\ \hline
		$2^{-12}$  &42  & 95.9 & 51 & 98.6 &33 &93.8& 58 & 92.0 & 61& 92.7&36 & 91.5 \\ \hline
		$2^{-11}$  & 42 & 95.9 & 54 & 98.4 &38 & 93.8& 60 & 91.1 & 59 &92.9 &33 &91.9 \\ \hline
		$2^{-10}$  & 42 & 96.2 & 54 & 97.3 &38 &94.0&65 & 91.0 &64 &92.6 & 36& 91.5  \\ \hline
		$2^{-9}$  &  42&  95.9& 56 & 97.5 &43 & 93.8&67 & 91.2 &62 &92.7 &36 & 91.4 \\ \hline
		$2^{-8}$  & 42 & 95.9 & 56 & 96.8 & 48&93.2 &60 & 91.9 &65 &92.5 &36&  91.4\\ \hline
		$2^{-7}$  & 42 & 95.9 & 59 & 96.8 &48 &91.3 &60 &91.7& 70& 92.1& 36& 91.3\\ \hline
		$2^{-6}$  & 42 & 95.9 & 59 & 96.8 &52 & 94.0& 65 & 92.2& 72&92.1 &38 & 91.6 \\ \hline
		$2^{-5}$  & 42 & 95.4 & 59 & 96.8 &52 &94.4 &58 &92.6& 74&92.2 &40 & 91.3 \\ \hline
		$2^{-4}$ & 42 & 95.9 & 59 & 97.3 & 57&93.8 &58&92.4&79 & 92.2&42 & 91.1\\ \hline
		$2^{-3}$  & 42 & 93.2 & 62 & 97.5 & 67&94.6 &63&91.1& 83&92.1 &40 & 91.7 \\ \hline
		$2^{-2}$  &50 & 89.7 & 62 & 97.7 & 67&94.4 &65&90.6&87 & 92.3&47 & 90.4 \\ \hline
		$2^{-1}$  & 50 & 92.7 & 64 & 98.2 & 71& 93.6& 67& 89.2 &90& 92.0& 51& 90.2 \\ \hline
		$2^{0}$ & 58 & 90.0 & 69 & 96.8 & 76&93.6 &71 & 88.1&91 & 91.9& 64& 87.6\\ \hline
		$2^{1}$ & 67 & 90.5 &77  & 95.2 & 81 &90.2& 75 &87.2  &92 &92.0 & 71& 85.4\\ \hline
		
		$2^{2}$  & 75 & 91.1 & 82 & 89.5 &81 & 88.5 &77 &82.6& 95& 91.8& 80& 80.8 \\ \hline
		$2^{3}$ & 83 & 88.6 & 90 & 76.4 &91 &73.6&83&77.3&100 &92.2 &91 &  68.2 \\ \hline
	\end{tabular}
\end{table}

\subsection{Results for global sparsity}
\label{Results for global sparsity}

This section is devoted to the global S-ORCT, i.e., when $\lambda^L=0$ and thus only global sparsity is taken into account. We focus on depth $D=2$, since for $D=1$ global sparsity is equal to local sparsity. Similarly to Subsection \ref{Results for local S-ORCTs}, Table \ref{Results2}

\begin{figure}
 \centering\hspace*{-1cm}
      \subfloat{
   \label{f:conejo}
    \includegraphics[width=0.36\textwidth]{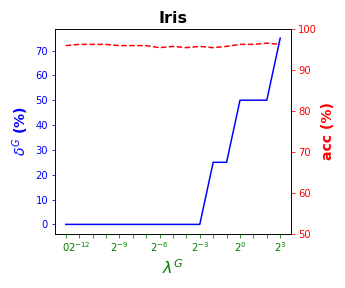}}
        \subfloat{
   \label{f:conejo}
    \includegraphics[width=0.36\textwidth]{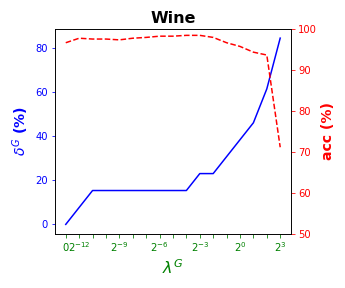}}
          \subfloat{
   \label{f:conejo}
    \includegraphics[width=0.36\textwidth]{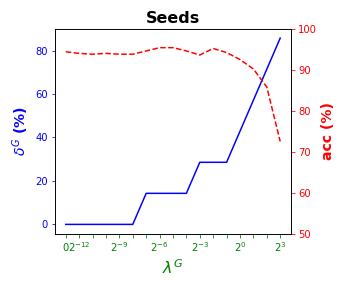}} \\
    \hspace*{-1cm}   \subfloat{
   \label{f:conejo}
    \includegraphics[width=0.36\textwidth]{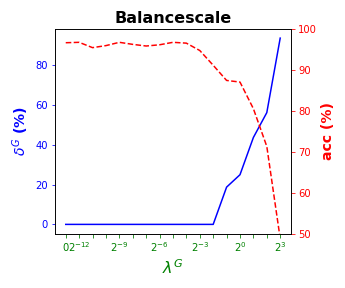}}
              \subfloat{
   \label{f:conejo}
    \includegraphics[width=0.36\textwidth]{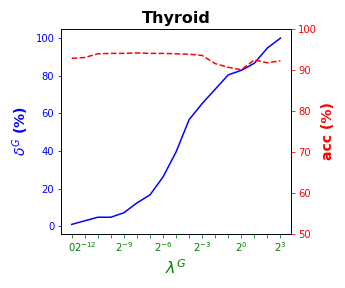}}
\subfloat{
	\label{f:conejo}
	\includegraphics[width=0.36\textwidth]{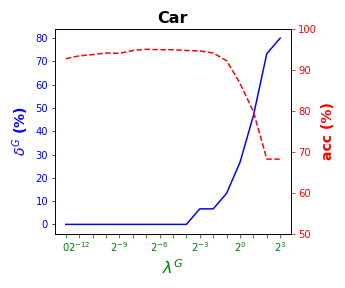}}
 \caption{ Graphical representation, for each data set, of the average percentage of predictor variables per tree, $\delta^G$, together with the average out-of-sample accuracy obtained, \textit{acc}, as a function of the values of $\lambda^G$ considered in the global S-ORCT construction.}
 \label{f:animales2}
\end{figure}

\begin{table}[h!]
	\caption{ Results for the global S-ORCT of depth $D=2$ as a function of $\lambda^G$, where $\delta^G$ represents the average percentage of predictor variables not used per tree over ten runs and \textit{acc}, the average out-of-sample accuracy.}
	\label{Results2}
	\centering
	\begin{tabular}{|c|c|c|c|c|c|c|c|c|c|c|c|c|}\hline
		\multirow{2}{*}{$\lambda^G$}  & \multicolumn{2}{|c|}{Iris} &\multicolumn{2}{|c|}{Wine} &\multicolumn{2}{|c|}{Seeds}& \multicolumn{2}{|c|}{ Balance}&\multicolumn{2}{|c|}{ Thyroid}&\multicolumn{2}{|c|}{Car} \\\cline{2-13}
		& $\delta^G$ & acc & $\delta^G$ & acc & $\delta^G$ & acc &  $\delta^G$ & acc & $\delta^G$ & acc  & $\delta^G$ & acc \\\hline
		$0$  & 0 & 95.9 & 0 & 96.6 & 0&94.4 & 0 & 96.6 & 1& 92.8& 0& 92.7 \\ \hline
		$2^{-12}$  & 0 & 96.2 & 18 & 97.7 &0 &94.0 &  0 &  96.7 & 3&93.0 &0 &  93.4  \\ \hline
		$2^{-11}$  &  0&  96.2& 15 & 97.5 & 0& 93.8 &  0 &  95.4 & 5& 93.9& 0&  93.7 \\ \hline
		$2^{-10}$  &  0&  96.2& 15 & 97.5 & 0&94.0 &  0 &  95.9 & 5&93.9 & 0& 94.1  \\ \hline
		$2^{-9}$  &  0&  95.9& 15 & 97.3 & 0&93.8 &  0 &  96.7&7 &94.0 &0 & 94.0 \\ \hline
		$2^{-8}$  &  0&  95.9& 15 &  97.7&0 & 93.8&  0 &  96.2& 12& 94.1&0 & 94.7 \\ \hline
		$2^{-7}$  &  0&  95.9& 15 & 97.9 &14 &94.6 &  0 &  95.8&17 &94.0 & 0&95.0 \\ \hline
		$2^{-6}$  &  0&  95.4& 15 & 98.2 & 14&95.4 &  0 &  96.1&26 &94.0 & 0&  94.9\\ \hline
		$2^{-5}$  &  2& 95.7 & 15 & 98.2 &14 &95.4 &  0 &  96.7&40 & 93.9& 0& 94.9 \\ \hline
		$2^{-4}$ &  0& 95.4 & 15 & 98.4  &14 &94.6 &  0 & 96.5& 57& 93.8&0 & 94.7\\ \hline
		$2^{-3}$  & 0 & 95.7 & 23 & 98.4 & 29& 93.6&  0 &  94.7 &65 & 93.5&7 &  94.6\\ \hline
		$2^{-2}$  &  25& 95.4 & 23 & 97.9 & 29& 95.2 &  0 &  91.1&73 & 91.5& 7& 94.1 \\ \hline
		$2^{-1}$  & 25 & 95.7 & 31 & 96.6 & 29&94.2 &  19 &  87.4&81 & 90.6&13 & 92.2 \\ \hline
		$2^{0}$ & 50 & 96.2 & 39 & 95.7 & 43& 92.5&  25 &  87.0& 83& 90.0& 27& 86.7\\ \hline
		$2^{1}$ & 50 & 96.2 & 46 & 94.3 & 57& 90.2& 44 &  80.5& 87& 92.4& 47& 79.8\\ \hline
		$2^{2}$  & 50 & 96.5 & 62 & 93.6 &71 &85.8 & 56 & 71.3&95 &91.7 & 73& 68.2 \\ \hline
		$2^{3}$ & 75 &  96.2& 85 & 71.1 & 86& 72.5&  94 &  48.8& 100& 92.2& 80&  68.2 \\ \hline
	\end{tabular}
\end{table}
 presents the results of the global S-ORCT, while Figure \ref{f:animales2}
  visualizes these results by showing simultaneously, per data set, $\delta^G$ (blue solid line) and \textit{acc} (red dashed line) as a function of the grid of the $\lambda^G$'s considered.
As for local sparsity, as $\lambda^G$ grows, $\delta^G$ increases. For Iris and Seeds, a similar classification accuracy to that with all of the predictor variables is obtained while removing the $75\%$ and $29\%$ of them, respectively. For Wine, the best rates of accuracy are obtained with $15\%-23\%$ of global sparsity. A loss of less than $10$ p.p. of accuracy is observed for Balance but $25\%$ of predictor variables are not being used, respectively. Car remains around the accuracy rate of $80\%$ while using half of the predictor variables. Thyroid, an imbalanced data set, is over the $90\%$ of accuracy for the whole grid of $\lambda^G$'s considered.

\subsection{Results for local and global sparsity}
\label{Results for local and global sparsity}

In this section, results enforcing local and global sparsity are presented by means of heatmaps, as seen in Figure \ref{heatmaps}. The experiment has been conducted on data sets of $K=3$ and $4$ classes, for which S-ORCTs of depth $D=2$ are built. For each dataset, three heatmaps are depicted as a function of the grid of the sparsity regularization parameters, $\lambda^L$ and $\lambda^G$: the average out-of-sample accuracy, \textit{acc}, and the local and global sparsities, $\delta^L$ and $\delta^G$, respectively, obtained over the ten runs performed. The color bar of each heatmap goes from light green to dark blue, being the latter the maximum accuracy, local sparsity or global sparsity achieved, respectively. As a general behavior, the best rates of accuracy are not always achieved only for $\left(\lambda^L,\lambda^G\right)=\left(0,0\right)$, but also for other pairs of the chosen grid, i.e., the data set remains equally well explained while needing less information. As before, according to local sparsity, for a fixed $\lambda^G$, $\delta^L$ has a growing trend. A similar behavior is observed for $\delta^G$ when $\lambda^L$ is fixed.  It is also worth mentioning that small changes of $\lambda^L$ quickly lead to a gain in $\delta^L$. Nevertheless, as expected, the gain in $\delta^G$ is slower for the same range in $\lambda^G$.

\begin{figure}
	\subfloat[Iris]{\hspace{-2cm}
	\includegraphics[scale=0.45]{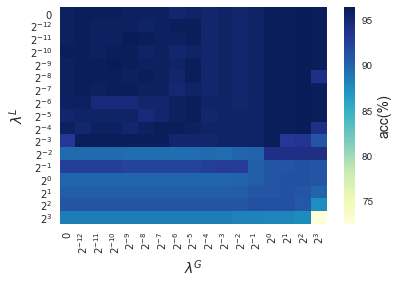}
	\includegraphics[scale=0.45]{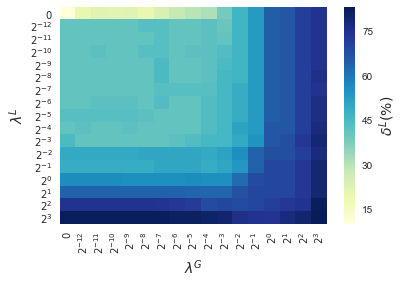}
	\includegraphics[scale=0.45]{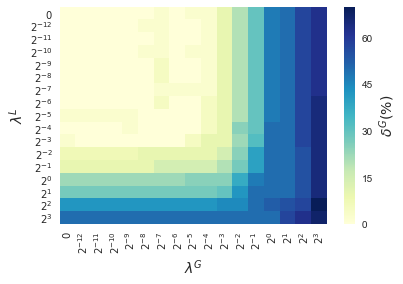}}\vspace{1.25cm}
	\subfloat[Wine]{\hspace{-2cm}
	\includegraphics[scale=0.45]{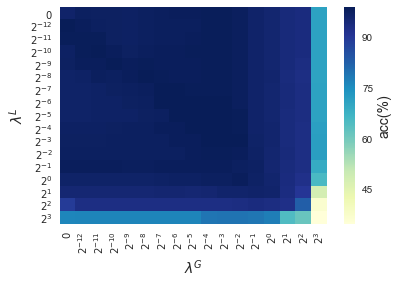}
	\includegraphics[scale=0.45]{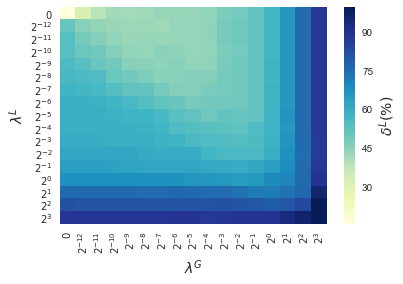}
	\includegraphics[scale=0.45]{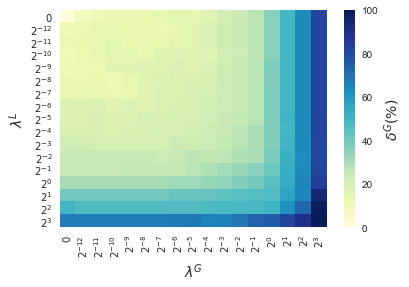}}\vspace{1.25cm}
	\subfloat[Seeds]{\hspace{-2cm}
	\includegraphics[scale=0.45]{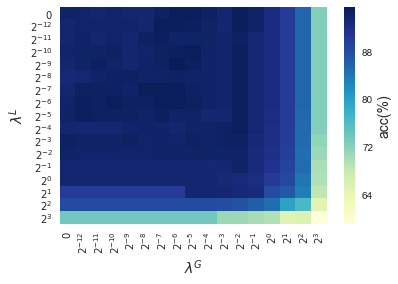}
	\includegraphics[scale=0.45]{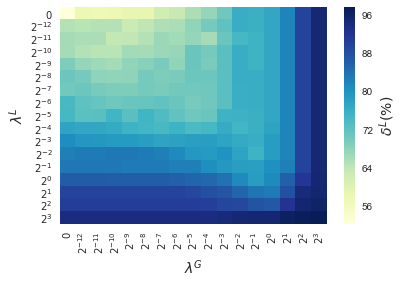}
	\includegraphics[scale=0.45]{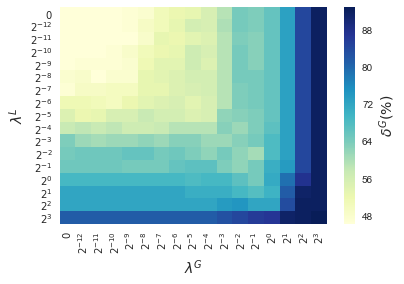}}
\phantomcaption
\end{figure}
\begin{figure}\ContinuedFloat
	\subfloat[Balance]{\hspace{-2cm}
	\includegraphics[scale=0.45]{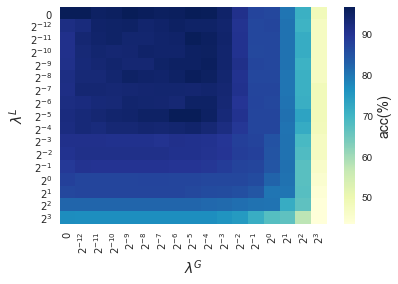}
	\includegraphics[scale=0.45]{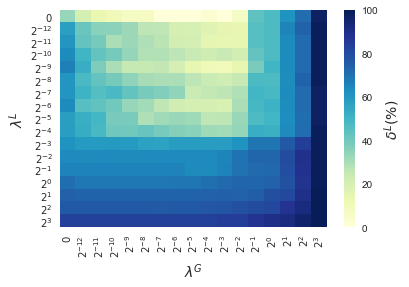}
	\includegraphics[scale=0.45]{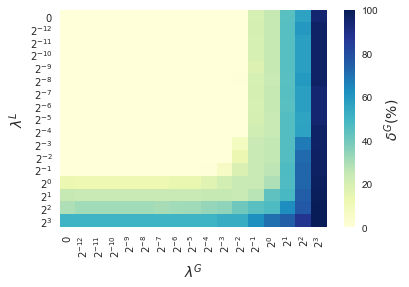}}\vspace{1.25cm}
\subfloat[Thyroid]{\hspace{-2cm}
	\includegraphics[scale=0.45]{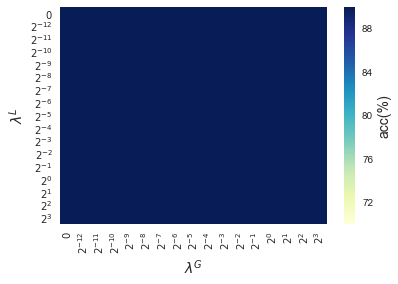}
	\includegraphics[scale=0.45]{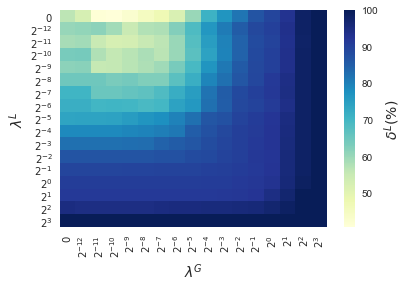}
	\includegraphics[scale=0.45]{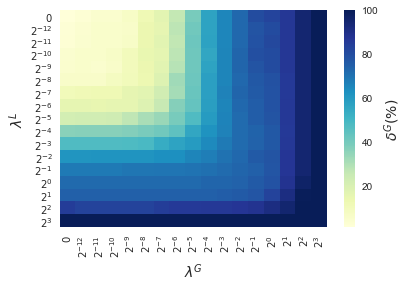}}\vspace{1.25cm}
\subfloat[Car]{\hspace{-2cm}
\includegraphics[scale=0.45]{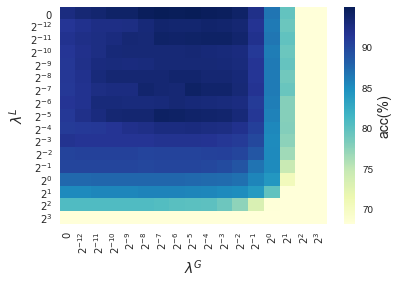}
\includegraphics[scale=0.45]{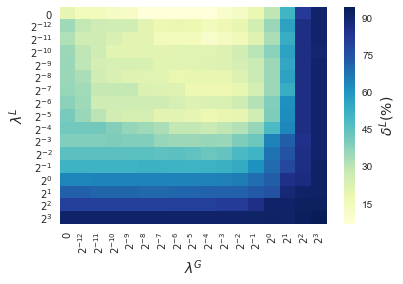}
\includegraphics[scale=0.45]{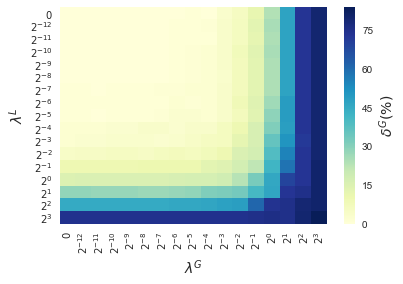}}\vspace{1.25cm}
\caption{ Heatmaps representation, for each data set, of the average out-of-sample accuracy, \textit{acc}, the average percentage of predictor variables not used per branch node, $\delta^L$, and the average percentage of predictor variables not used per tree, $\delta^G$, respectively, as a funcion of the grid of the sparsity parameters, $\lambda^L$ and $\lambda^G$, considered in the S-ORCT of depth $D=2$ construction.
}
\label{heatmaps}
\end{figure}

\subsection{Comparison S-ORCT versus CART}
\label{Comparison S-ORCT versus CART}

A statistical comparison between the proposed S-ORCT and CART, the classic approach to build decision trees, is provided in this section. As stated in the introduction of the paper, CARTs, as many other approaches that implement orthogonal cuts \cite{bertsimas2017optimal,firat2018constructing,menickelly2016optimal}, are leaders in terms of local sparsity. Thus, the comparison S-ORCT versus CART is performed in terms of accuracy and global sparsity. Tables \ref{Results for $D=1$.} and \ref{Results2} for S-ORCT have been considered for the experiment.

CART has been trained and tested over the same ten runs as S-ORCT. For each pair S-ORCT$\left(\lambda^G\right)$ versus CART, two hypothesis tests for the equality of means of paired samples were carried out, one for accuracy and another for global sparsity, assuming normality, at a $5\%$ significance level. For this task, the \texttt{t.test} function in R has been used. Figure \ref{statisticalsignificance2} depicts, for each data set, the resulting confidence intervals (blue solid line) at the $95\%$ confidence level for the difference in average accuracy (on the left) and global sparsity (on the right) between S-ORCT$\left(\lambda^G\right)$ and CART. The red dashed horizontal line represents the null hypothesis in each case. Except for Creditapproval and Thyroid, for the smaller values of $\lambda^G$, our approach is significantly better than, or at least comparable to, CART in terms of accuracy, while CART is significantly better than, or at least comparable to, in terms of global sparsity. For the larger values of $\lambda^G$, our approach starts to be comparable and then dominate CART in terms of global sparsity at the cost of accuracy.

\begin{figure}\vspace{-1cm} 	
	\centering
	\subfloat[Monks-3]{
		\includegraphics[scale=0.8]{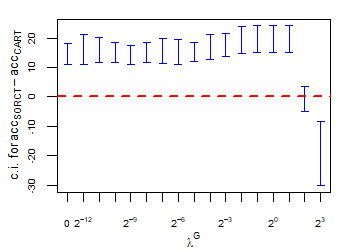}\hspace*{1cm}
		\includegraphics[scale=0.8]{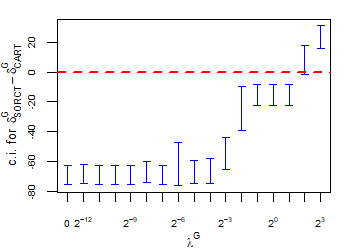}}\\
	\subfloat[Monks-1]{
		\includegraphics[scale=0.8]{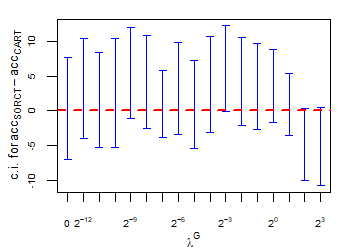}\hspace*{1cm}
		\includegraphics[scale=0.8]{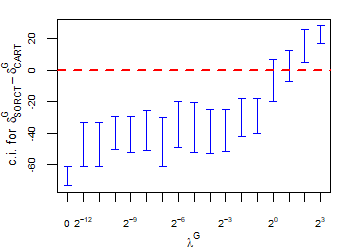}}\\
	\subfloat[Monks-2]{
		\includegraphics[scale=0.8]{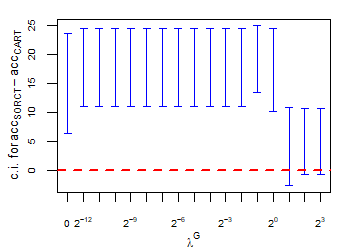}\hspace*{1cm}
		\includegraphics[scale=0.8]{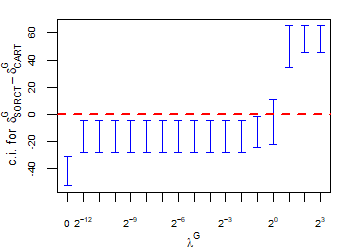}}\\
	\subfloat[Sonar]{
		\includegraphics[scale=0.8]{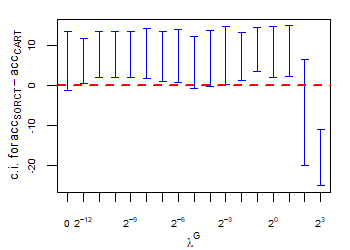}\hspace*{1cm}
		\includegraphics[scale=0.8]{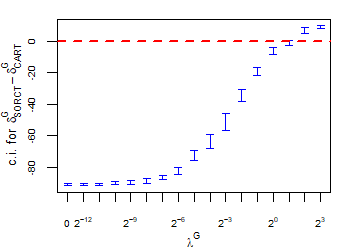}}
	\phantomcaption
\end{figure}
\begin{figure}\ContinuedFloat\vspace{-1cm} 	
	\centering
	\subfloat[Ionosphere]{
		\includegraphics[scale=0.8]{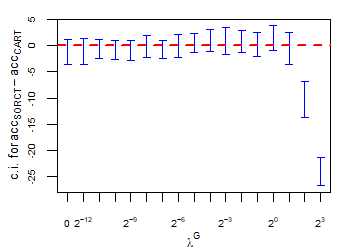}\hspace*{1cm}
		\includegraphics[scale=0.8]{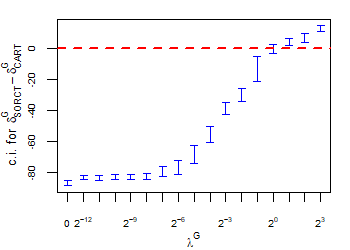}}\\
	\subfloat[Wisconsin]{
		\includegraphics[scale=0.8]{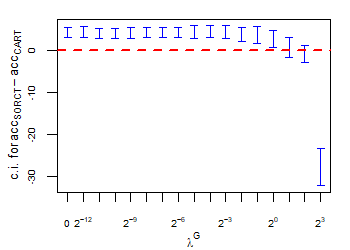}\hspace*{1cm}
		\includegraphics[scale=0.8]{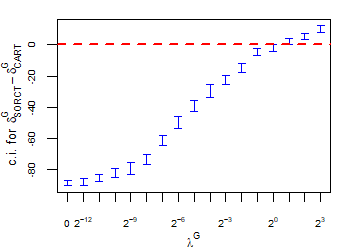}}\\
	\subfloat[Creditapproval]{
		\includegraphics[scale=0.8]{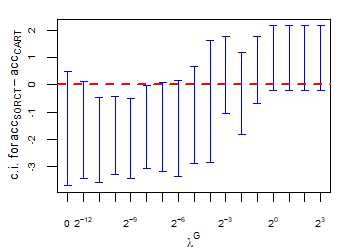}\hspace*{1cm}
		\includegraphics[scale=0.8]{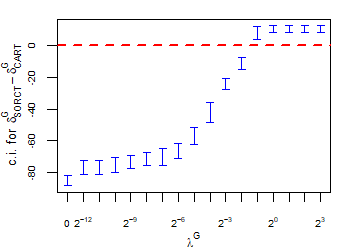}}\\
	\subfloat[Pima]{
		\includegraphics[scale=0.8]{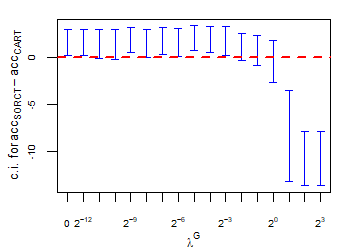}\hspace*{1cm}
		\includegraphics[scale=0.8]{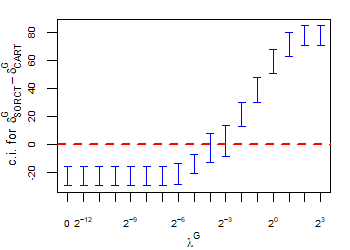}}
	\phantomcaption
\end{figure}
\begin{figure}\ContinuedFloat\vspace{-1cm} 	
	\centering
	\subfloat[Germancredit]{
		\includegraphics[scale=0.8]{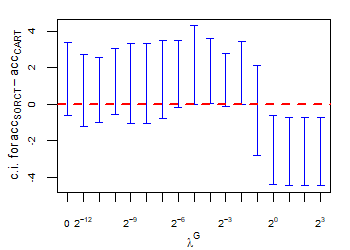}\hspace*{1cm}
		\includegraphics[scale=0.8]{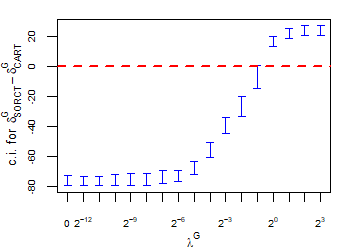}}\\
	\subfloat[Banknote]{
		\includegraphics[scale=0.8]{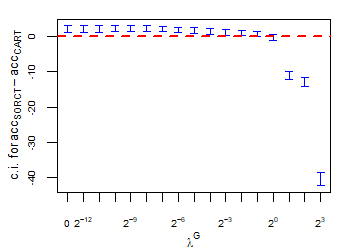}\hspace*{1cm}
		\includegraphics[scale=0.8]{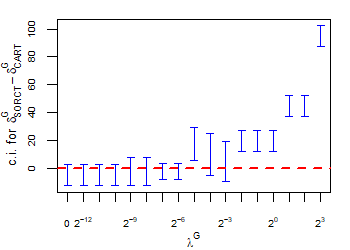}}\\
	\subfloat[Ozone]{
		\includegraphics[scale=0.8]{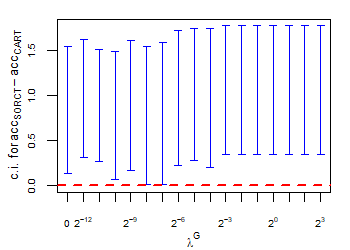}\hspace*{1cm}
		\includegraphics[scale=0.8]{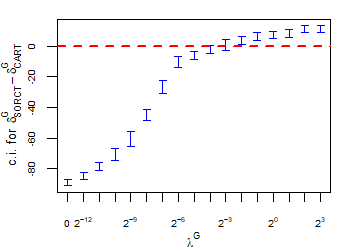}}\\
	\subfloat[Spam]{
		\includegraphics[scale=0.8]{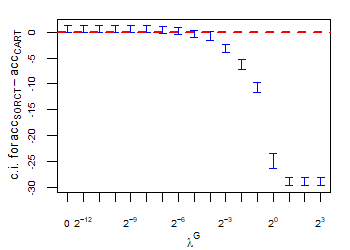}\hspace*{1cm}
		\includegraphics[scale=0.8]{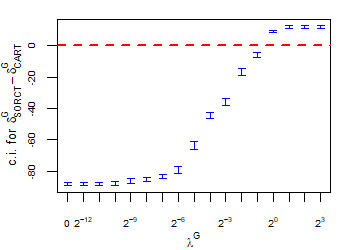}}
	\phantomcaption
\end{figure}
\begin{figure}\ContinuedFloat\vspace{-1cm} 	
	\centering
	\subfloat[Iris]{
		\includegraphics[scale=0.8]{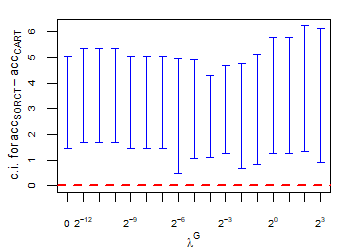}\hspace*{1cm}
		\includegraphics[scale=0.8]{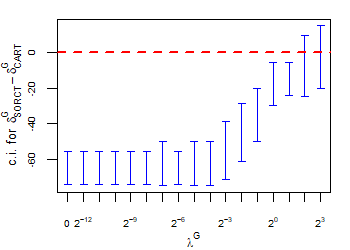}}\\
	\subfloat[Wine]{
		\includegraphics[scale=0.8]{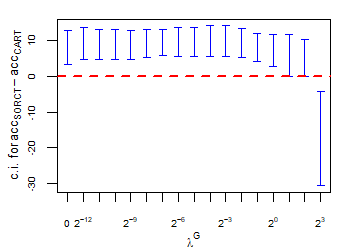}\hspace*{1cm}
		\includegraphics[scale=0.8]{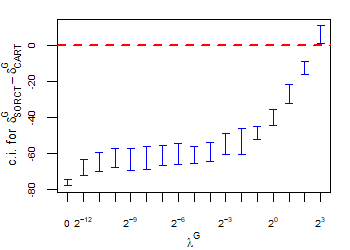}}\\
	\subfloat[Seeds]{
		\includegraphics[scale=0.8]{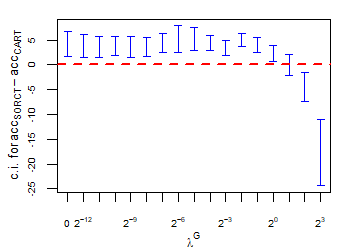}\hspace*{1cm}
		\includegraphics[scale=0.8]{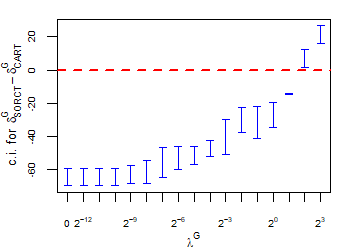}}\\
	\subfloat[Balance]{
		\includegraphics[scale=0.8]{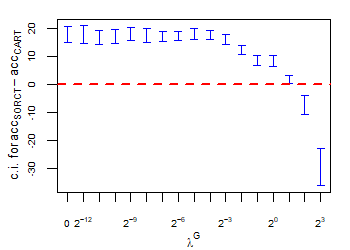}\hspace*{1cm}
		\includegraphics[scale=0.8]{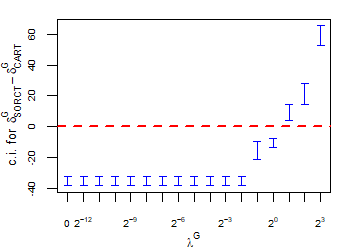}}
	\phantomcaption
\end{figure}
\begin{figure}\ContinuedFloat\vspace{-1cm} 	
	\centering
	\subfloat[Thyroid]{
		\includegraphics[scale=0.8]{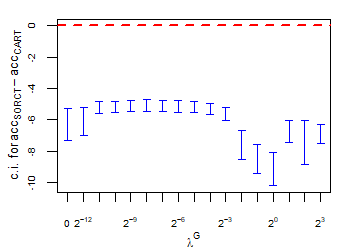}\hspace*{1cm}
		\includegraphics[scale=0.8]{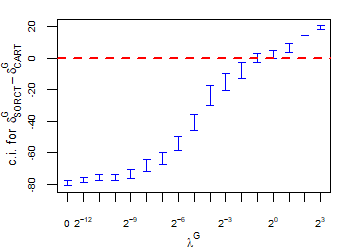}}\\
	\subfloat[Car]{
		\includegraphics[scale=0.8]{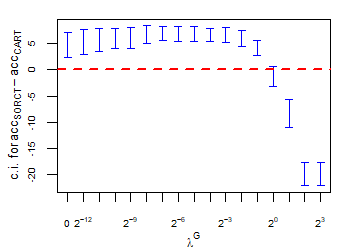}\hspace*{1cm}
		\includegraphics[scale=0.8]{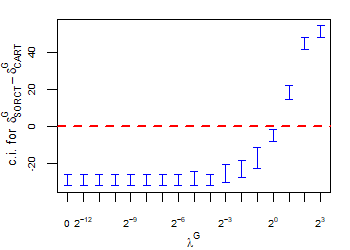}}
\caption{ Graphical representation, for each data set, of the confidence intervals (blue solid line) at the $95\%$ for the difference in average accuracy (on the left) and global sparsity (on the right) between S-ORCT$\left(\lambda^G\right)$ and CART. The red dashed horizontal line represents the null hypothesis in each case.} \label{statisticalsignificance2}
\end{figure}

\section{Conclusions and future research}
\label{Conclusions and future research}

Recently, several proposals focused on building optimal classification trees are found in the literature to address the shortcomings of the classic greedy approaches. In this paper, we have proposed a novel continuous optimization-based approach, the Sparse Optimal Randomized Classification Tree (S-ORCT), in which a compromise between good classification accuracy and sparsity is pursued. Local and global sparsity in the tree are modeled by including in the objective function norm-like regularizations, namely, $\ell_1$ and $\ell_\infty$, respectively. Our numerical results illustrate that our approach can improve both sparsities without harming classification accuracy. Unlike CART, we are able to easily trade in some of our classification accuracy for a gain in global sparsity.

Some extensions of our approach are of interest. First, this metholodogy can be extended straightaway to a regression tree counterpart, where the response variable is continuous. Second, categorical data is addressed in this paper through the inclusion of dummy predictor variables. For a given categorical predictor variable, and by means of an $\ell_\infty$-norm regularization, one can link all its dummies across all the branch nodes in the tree, with the aim of better modeling its contribution to the classifier. Third, it is known that bagging trees tends to enhance accuracy. An appropiate bagging scheme of our approach, where sparsity is a key point, is a nontrivial design question.

\par\noindent  {\bf Acknowledgements.} This research has been financed in part by research projects EC H2020 MSCA RISE NeEDS (Grant agreement ID: 822214), COSECLA - Fundaci\'on BBVA, MTM2015-65915R, Spain, P11-FQM-7603 and FQM-329, Junta de Andaluc\'{\i}a, the last three with EU ERF funds.  This support is gratefully acknowledged.


\bibliographystyle{apa}
\bibliography{elsarticle-template_revision2}

\end{document}